\documentclass{article}

\usepackage{arxiv}

\usepackage[utf8]{inputenc} 
\usepackage[T1]{fontenc}    
\usepackage{hyperref}       
\usepackage{url}            
\usepackage{booktabs}       
\usepackage{amsfonts}       
\usepackage{nicefrac}       
\usepackage{microtype}      
\usepackage{lipsum}		    
\usepackage{graphicx}
\usepackage[square,sort,comma,numbers]{natbib}
\usepackage{doi}
\usepackage{float}

\usepackage{caption}
\usepackage{subcaption}
\usepackage{xcolor}
\usepackage{amsmath}

\renewcommand{\i}{\text{i}} 
\newcommand{\R}{{\rm I\!R}} 
\newcommand{\F}[1]{\mathcal{F}\left \{#1 \right\}} 
\newcommand{\IF}[1]{\mathcal{F}^{-1} \left \{#1 \right \}} 
\renewcommand{\O}[1]{\mathcal{O} (#1 )} 

\title{Solving the complete pseudo-impulsive radiation and diffraction problem using a spectral element method}

\author{ \href{https://orcid.org/0000-0001-6698-2623}{\includegraphics[scale=0.06]{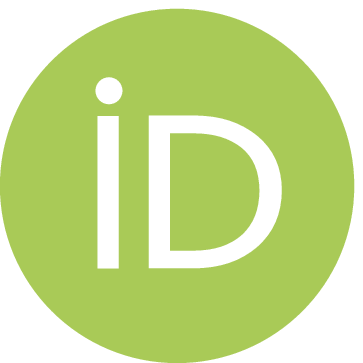}\hspace{1mm}Jens Visbech} \\
	Department of Applied Mathematics and Computer Science\\
	Technical University of Denmark\\
	Kongens Lyngby, 2800 \\
	\texttt{jvis@dtu.dk} \\
	\And
	\href{https://orcid.org/0000-0001-8626-1575}{\includegraphics[scale=0.06]{Figures/orcid.eps}\hspace{1mm}Allan P. Engsig-Karup} \\
	Department of Applied Mathematics and Computer Science\\
	Technical University of Denmark\\
	Kongens Lyngby, 2800 \\
	\texttt{apek@dtu.dk} \\
	\And
	\href{https://orcid.org/0000-0002-7263-442X}{\includegraphics[scale=0.06]{Figures/orcid.eps}\hspace{1mm}Harry B. Bingham} \\
	Department of Civil \& Mechanical Engineering\\
	Technical University of Denmark\\
	Kongens Lyngby, 2800 \\
	\texttt{hbbi@dtu.dk} \\
}

\renewcommand{\shorttitle}{An \textit{arXiv} preprint}

\hypersetup{
pdftitle={},
pdfsubject={},
pdfauthor={Jens Visbech, Allan Peter Engsig-Karup, Harry Bingham},
pdfkeywords={},
}

\begin{document}
\maketitle

\begin{abstract}
	This paper presents a novel, efficient, high-order accurate, and stable spectral element-based model for computing the complete three-dimensional linear radiation and diffraction problem for floating offshore structures. We present a solution to a pseudo-impulsive formulation in the time domain, where the frequency-dependent quantities, such as added mass, radiation damping, and wave excitation force for arbitrary heading angle, $\beta$, are evaluated using Fourier transforms from the tailored time-domain responses. The spatial domain is tessellated by an unstructured high-order \textit{hybrid} configured mesh and represented by piece-wise polynomial basis functions in the spectral element space. Fourth-order accurate time integration is employed through an explicit four-stage Runge-Kutta method and complemented by fourth-order finite difference approximations for time differentiation. To reduce the computational burden, the model can make use of symmetry boundaries in the domain representation. The key piece of the numerical model -- the discrete Laplace solver -- is validated through $p$- and $h$-convergence studies. Moreover, to highlight the capabilities of the proposed model, we present \textit{prof-of-concept} examples of simple floating bodies (a sphere and a box). Lastly, a much more involved case is performed of an oscillating water column, including generalized modes resembling the piston motion and wave sloshing effects inside the wave energy converter chamber. In this case, the spectral element model trivially computes the infinite-frequency added mass, which is a singular problem for conventional boundary element type solvers.
\end{abstract}

\keywords{Spectral element method \and High-order numerical method \and Free surface \and Potential flow \and Radiation \and Diffraction \and Floating offshore structures \and Oscillating water column \and Hybrid mesh \and Pseudo-impulse}

\section{Introduction}\label{sec:introdcution}

The offshore industry encompasses a broad spectrum of energy-producing activities, particularly those focusing on green and renewable sources. Wind, solar, and wave energy are sustainable energy resources that can be harnessed  from the ocean environment. This is done via floating offshore structures, which have traditionally been developed and optimized in large laboratory wave flumes. However, physical exercises of this nature can be time-consuming and expensive. To simplify the process, numerical models have become increasingly popular over the past decades due to advances in computational resources, e.g., through increasingly efficient algorithms and methods implemented on high-performance computing systems. The dynamics are modeled through a set of \textit{governing equations} that mimic physical laws, which are then discretized using \textit{numerical methods}. This introduction focuses solely on the numerical modeling of ocean waves and their interactions with offshore structures. For a recent comprehensive modeling review for wave energy converters, please consider \cite{davidson2020efficient}.

\subsection{On modeling wave propagation and wave-structure interactions}
As for the governing equations, motions in fluid flows are known to be well-described by the Navier-Stokes equations (NSE). However, using such high-fidelity equations comes at a cost in terms of computational expense due to the complex nonlinear structure of the mass and momentum formulations. Applying physical constraints to the NSE can dramatically reduce this expense, creating a new set of governing equations to model the fluid motion. Assuming that the flow is of constant density (incompressible), without effects from viscosity (inviscid), and free of a net rotational motion (irrotational), the NSE can be reduced to the potential flow (PF) equations \cite{engsigkarup2013fast}. Here, an elliptic problem governs the spatial domain -- in terms of the Laplace equation -- combined with a coupled hyperbolic system that controls the temporal domain's evolution. Both nonlinear and linear versions of the PF formulations are feasible. Due to the reduced computational requirements, such formulations are well-suited for large-scale and long-time simulations.

The linear version is obtained by assuming small amplitude waves (compared to the wavelength) and small structural movements (concerning the length scale of the floating body). With this, it is sufficient to represent the wave dynamics solely using first-order terms from the nonlinear perturbation series, \cite{svendsen2001hydrodynamics}. The linear PF formulation can be described in the frequency domain, where a typical task is to determine some arbitrary system dynamic due to a single wave frequency. However, using a time domain formulation, this computational exercise can solve for a wider wave frequency span in one simulation. One approach is called the impulsive formulation and was studied extensively in the 1980s and 1990s, e.g., consider \cite{liapis1986time,beck1989time,beck1994time,korsmeyer1998forward}, based on the work done by \cite{cummins1962impulse} in the 1960s. The key idea is to impose a rapid initial velocity impulse of the body in the time domain and let the system evolve. Hereafter, the frequency content can be obtained via Fourier transformations. Yet, as initially noticed in \cite{king1987time} and numerically shown in \cite{visbech2023spectral}, the impulsive formulation will impose unit amplitude energy at all frequencies independent of the spatial discretization. Ultimately, this leads to the phenomenon of spurious oscillations in the system. A solution to this problem is exploited by the pseudo-impulsive formulation, \cite{king1987time}, where the discretization governs the frequency content by the use of a tailored Gaussian displacement function.

\subsection{On numerical methods for potential flow formulations} 
Considering the PF formulation (both linear and nonlinear) for modeling wave propagation and wave-structure interactions, various numerical methods have been employed to solve those equations discretely. The ideal  scenario is to have a method that is: i) sufficiently accurate in time and space, ii) stable with respect to time integration, iii) numerically efficient and scalable, and iv) able to handle complex geometrical tasks. Below, we briefly outline a few major developments in numerical methods related to PF modeling.

\textit{The boundary element method} (BEM): Through Green's 2nd identity, the BEM recasts the Laplace equation into a boundary integral equation, where the volumetric information is projected onto the surface of the boundary of interest. With this, the dimension of the problem is reduced by one; however, at the cost of a loss of sparsity. The BEM employs a discrete panel/element-based approach to approximate geometrical features, making the method ideal for problems with complex bodies and bathymetry. The computational bottleneck of the problem is solving the dense linear system of equations associated with the Laplace problem. This task can naively be performed with $\mathcal{O}(n^3)$ -- with $n$ being degrees of freedom -- in computational scaling using direct solvers and $\mathcal{O}(n^2)$ using iterative solvers. This (critical) scaling issue can be improved using a fast multipole approach due to \cite{greengard1987fast} to $\mathcal{O}(n \log (n))$; however, with a large constant in front. Further improvements (more than one order of magnitude at engineering accuracy) can be made using a pre-corrected fast Fourier transform (FFT) method due to \cite{phillips1997precorrected}. For PF problems, many works exist describing the use of BEM, e.g., the historical review see \cite{cheng2005heritage} and recent works \cite{harris2022nonlinear,yan2011efficient}. Commercially, the \textit{state-of-the-art} tool WAMIT employs the BEM when simulating linear and nonlinear wave-structure interactions in the frequency domain, \cite{lee2006wamit}. Open source options are also available, e.g., through NEMOH, \cite{babarit2015theoretical,kurnia2023nemoh}, that solves first- and second-order hydrodynamic loads in the frequency domain. Impulsive time domain formulations are presented in \cite{kara2011time}.

\textit{The finite difference method} (FDM): Contrary to the BEM, the FDM solves the Laplace problem on the entirety of the computational domain using locally applied stencils. This yields great possibilities for the use of low- and high-order accurate schemes with efficient linear scaling. The incorporation of very complex geometrical features sets requirements for the method, which does not naturally follow.  A low- and high-order accurate FDM was applied for pure nonlinear wave propagation in \cite{bingham2007accuracy,engsigkarup2009efficient,li1997three} using $\sigma$-transformations to cope with varying bathymetry. Also, nonlinear wave-structure interactions have been studied in \cite{ducrozet2010high}. Open source options are available through OceanWave3D \cite{engsigkarup2009efficient}. In relation to linear pseudo-impulsive formulations, some work has been done in the 2010s for the radiation and diffraction problem with and without zero-speed, \cite{read2012overset,amini2017solving,amini2018pseudo}.

\textit{The finite element method} (FEM): Within engineering (especially mechanical engineering), the FEM is one of the most used numerical discretization methods, \cite{zienkiewicz2005finite}. The Galerkin-based method solves the governing equations -- in a weak sense -- on a discrete domain. The domain is tessellated by finite elements that take various shapes and have different properties, making the method extremely capable of solving complex geometrical problems. The FEM uses piece-wise linear functions as the expansion basis with local support, resulting in a sparse linear system of equations that converge with second-order accuracy. This allows for efficient linear-scaled, $\mathcal{O}(n)$, iterative solution strategies. For nonlinear wave propagation and structure interactions, see \cite{wu1994finite,taylor1996analysis,wu2003coupled}.

\textit{The spectral element method} (SEM): The SEM is a natural extension of the FEM, which substitutes the linear basis function with polynomials of arbitrary order, hereby enabling the accuracy of the method to match -- or even exceed -- the spatial dimension of the PF problem, where low-order methods fall short according to \cite{engsigkarup2017stabilsedpart2}. The idea of the SEM was first introduced by \cite{patera1984spectral} for fluid flows. The SEM employs the finite element-based ability to converge by refining the mesh ($h$-convergence) but also by increasing the order of the expansion basis ($p$-convergence). With this, yielding a very efficient scheme -- $\mathcal{O}(n)$ with iterative solvers as shown in \cite{engsigkarup2021efficient} -- that naturally copes with complex geometries. In the context of free surface flows, the SEM was shown to be inherently unstable by \cite{robertson1999free} in the late 1990s; however, almost two decades later, the first stable nonlinear model for wave propagation was presented by \cite{engsigkarup2016stabilised}. Since then, considerable work has been put into modeling wave-wave and wave-structure interactions. To name a few, see \cite{engsigkarup2018spectral,mortensen2021simulation,engsigkarup2019mixed,bosi2019spectral} and consider a recent review of SEM modeling by \cite{xu2018spectral}.

\subsection{Paper contributions}

This work presents a novel SEM-based model for solving the linear three-dimensional radiation and diffraction problem with zero speed through a pseudo-impulsive formulation. The entire model is high-order accurate as i) the spatial discretization is treated with piece-wise polynomial basis functions of arbitrary order, and ii) the temporal domain employs a fourth-order accurate Runge-Kutta time integration strategy and fourth-order accurate finite difference stencils for time differentiation. Moreover, hydrodynamic quantities are mapped efficiently between the time and the frequency domain by the use of FFTs and their inverse, respectively. The discrete high-order model employs so-called \textit{hybrid} meshes, allowing for stable time integration and modeling complex floating offshore problems. Moreover, the model has the possibility to employ up to two symmetry boundaries to reduce the computational cost. The paper considers $p$- and $h$-convergence of the governing Laplace operator. Furthermore, a \textit{proof-of-concept} validation study of a floating sphere and a floating box is carried out for the complete velocity potential to establish the legitimacy of the proposed solution strategy. Lastly, a more involved case study is carried out of an oscillating water column (OWC) with the use of generalized modes on a special internal boundary of the chamber. 

\subsubsection{Outline}

The remainder of the paper is structured in the following way: In Section \ref{sec:mathematical_problem}, the mathematical problem is stated by presenting the computational domain and the governing equations for the entire linear pseudo-impulsive formulation. Hereafter, in Section \ref{sec:discretization_methods} outlines the numerical solution approach in terms of numerical methods and discretization. Section \ref{sec:validation_and_computational_properties} validates the numerical model. Finally, in Section \ref{sec:numerical_results}, various numerical results are presented of a floating sphere, a box, and a full-scale floating wave energy converter in terms of the OWC. At last, the key points are summarized in Section \ref{sec:conclusion}. The Appendix is located in Section \ref{sec:appendix}.
\section{Mathematical Problem}\label{sec:mathematical_problem}

\subsection{Preliminaries - The computation domains and definitions}

The fluid domain is of three spatial dimensions ($d = 3$), defined in Cartesian coordinates spanned by $\boldsymbol{x} = (x,y,z)$, and denoted by $\Omega \in \R^d$. The combined $x$- and $y$-axis covers the horizontal plane, whereas the $z$-axis covers the finite depth in the vertical direction defined positively upwards. The volumetric domain is bounded completely in a Lipschitz sense by $\partial \Omega =\Gamma \in \R^d$. The boundary is partitioned into: the undisturbed free surface at $z=0$, $\Gamma^{\text{FS}}$, the bathymetry, $\Gamma^{\text{b}}$, the far-field, $\Gamma^{\infty}$, the body, $\Gamma^{\text{body}}$, and lastly -- if the geometrical features of $\Omega$ allow for it -- symmetry boundaries in both horizontal directions, $\Gamma^{\text{sym},x}$ and $\Gamma^{\text{sym},y}$. The former of the symmetry boundaries is defined about the $x$-axis (the $y=0$ plane), whereas the latter is defined about the $y$-axis (the $x = 0$ plane). We note the possibility of more -- so-called special -- boundaries, $\bigcup_{i = 1}^n \Gamma^{\text{s},i}$ than those aforementioned six, which will be highlighted in Section \ref{sec:OWC}. With this, $\Gamma = \Gamma^{\text{FS}} \cup \Gamma^{\text{b}} \cup \Gamma^{\infty} \cup \Gamma^{\text{body}} \cup \Gamma^{\text{sym,x}} \cup \Gamma^{\text{sym,y}} \cup \left \{ \bigcup_{i = 1}^n \Gamma^{\text{s},i} \right \}$.

\begin{figure}[t] 
    \centering
    \includegraphics{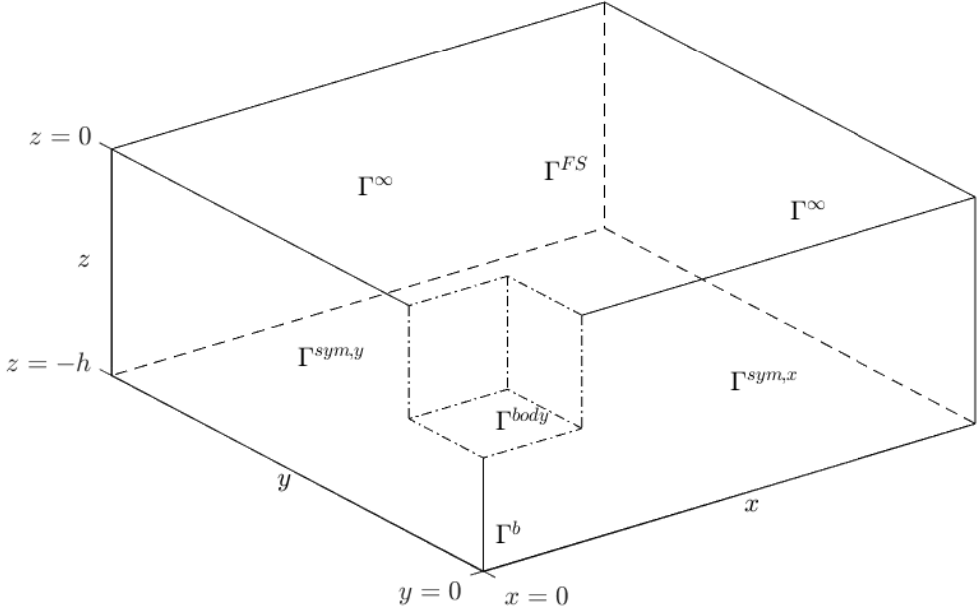}
    \caption{Concept of the three-dimensional fluid domain.}
    \label{fig:fluid_domain}
\end{figure}

The temporal domain is denoted by $\mathcal{T} : t \geq 0$ and its frequency-based counterpart by $\mathcal{H} : \omega \geq 0$, where $\omega = \frac{f}{2 \pi}$ is the radian frequency connected to the cyclic frequency, $f$. Concerning $\Omega$, we define some physical quantities, namely the free surface elevation, $ \eta = \eta(x,y,t): \Gamma^{\text{FS}} \times \mathcal{T} \mapsto \R$ and the depth of the water (measured from $z=0$ to the seabed), $h = h(x,y): \Gamma^b \mapsto \R$; however, for the entirety of this work, we take $h$ to be constant. See Figure \ref{fig:fluid_domain} for a conceptual layout of $\Omega$ and its boundaries.

To model the fluid interactions between ocean waves and floating offshore structures, we adopt the setting of potential flow theory, where the fluid is considered to be: i) incompressible, ii) inviscid, and iii) irrotational. With this, the fluid velocity field can be expressed in terms of the gradient of a scalar velocity potential, $\boldsymbol{u} = \nabla \phi$, where $\boldsymbol{u} = (u,v,w) =  (u(\boldsymbol{x},t),v(\boldsymbol{x},t),w(\boldsymbol{x},t)): \Omega \times \mathcal{T} \mapsto \R^3$ and $\phi = \phi(\boldsymbol{x},t) : \Omega \times \mathcal{T} \mapsto \R$ with $\nabla = \left (\partial_x, \partial_y, \partial_z \right )$ being the Cartesian gradient operator in all three spatial directions. We arrive at the linearized potential flow setting by applying further simplification to the potential flow formulation, namely only accounting for small amplitude waves and structure movements by discarding all but first-order terms in the nonlinear free surface perturbation. Within this setting, we can decompose the complete velocity potential into independently solvable parts due to \cite{haskind1946hydrodynamic} as
\begin{equation}
    \phi = \phi_0 + \phi_s + \sum_{k=1}^{6+K} \phi_k
\end{equation}
where $\phi_0$ and $\phi_s$ constitute the diffraction potential. The diffraction potential is composed of the known incident potential (former) and the unknown scattered potential (latter) due to fixed body interactions with the incident waves. Moreover, $\phi_k$ for $1 \leq k \leq 6$, are the six radiation potentials due to rigid-body movements in three translatory degrees of freedom (surge: $k = 1$, sway: $k = 2$, and heave: $k = 3$) in the $x$-, $y$-, and $z$-direction, respectively, and three rotational degrees of freedom (roll: $k = 4$, pitch: $k = 5$, and yaw: $k = 6$) around the $x$-, $y$-, and $z$-axis, respectively. For $k > 6$, $\phi_k$ denotes generalized modes of motion, as will become apparent in Section \ref{sec:OWC} in connection with the aforementioned special boundaries. Until further notice, $K = 0$.

\subsubsection{Other definitions}

We define the Fourier transform operator as $\F{\cdot}: \mathcal{T} \mapsto \mathcal{H}$ and its inverse as $\IF{\cdot}: \mathcal{H} \mapsto \mathcal{T}$ with respect to the angular frequency, $\omega$. For any periodic function, $g(t)$, we have that
\begin{align}\label{eq:FT}
\hat{g}(\omega) &= \F{g(t)} = \int_{-\infty}^{\infty} g(t) e^{-i \omega t} dt, \quad \text{in} \quad \mathcal{H}, \\
g(t) &= \IF{\hat{g}(\omega)} = \frac{1}{2\pi} \int_{-\infty}^{\infty} \hat{g}(\omega) e^{i \omega t} d\omega, \quad \text{in} \quad \mathcal{T}.\label{eq:IFT}
\end{align}

\subsection{Linear free surface potential flows}\label{sec:linear_potential_flow}

The following outlines the complete mathematical problem and the governing equations for the linear pseudo-impulsive radiation and diffraction problems. First, in Section \ref{sec:linear_potential_flow}, the more universal equations for linear free surface potential flows are presented, whereas special conditions and treatments for the radiation and diffraction problems are outlined separately in Section \ref{sec:radiation} and \ref{sec:diffraction}, respectively.

Now, considering any decomposed part, $\phi_i$, of the complete velocity potential, $\phi$, where $i = \{s,0,1,...,K\}$. The spatial governing scalar Laplace equation is at the heart of any potential flow formulation. The computation task is to solve this mass-conserving problem such that $\phi_i \in C^2(\Omega)$
\begin{equation}
     \nabla^2 \phi_i  = 0, \quad \text{in} \quad \Omega,
\end{equation}
with boundary conditions as
\begin{align}
     \phi_i &= \phi_i, \quad \text{on} \quad \Gamma^{\text{FS}}, \\
     \partial_n \phi_i &= 0, \quad \text{on} \quad \Gamma^{\text{b}} \quad \text{and} \quad \Gamma^{\infty}, \\
     (1-\theta) \phi_i + \theta \partial_n \phi_i &= 0 , \quad \text{on} \quad \Gamma^{\text{sym,x}}, \quad \text{if} \quad \Gamma^{\text{sym,x}} \notin \emptyset, \\
     (1-\theta) \phi_i + \theta \partial_n \phi_i &= 0 , \quad \text{on} \quad \Gamma^{\text{sym,y}},  \quad \text{if} \quad \Gamma^{\text{sym,y}} \notin \emptyset,
\end{align}
where $\partial_n$ is the derivative with respect to the normal direction of the boundary. Also, on symmetry boundaries, $\Gamma^{\text{sym,x}}$ and $\Gamma^{\text{sym,y}}$, one can impose homogeneous first-type Dirichlet boundary conditions (if the solution is anti-symmetric about this boundary) or homogeneous second-type Neumann boundary conditions (if the solution is symmetric about this boundary) by setting the Boolean parameter $\theta$. Specific $\theta$ values are specified in Section \ref{sec:radiation} for the radiation problems and in Section \ref{sec:diffraction} for the diffraction problem. To make the boundary value problem well-posed, conditions on $\Gamma^{\text{body}}$ are required. Those conditions really define and differentiate the radiation from the diffraction problem, thus outlined in the two separate sections.

The boundary value problem connects to the coupled linearized Eulerian time-dependent dynamic and kinematic free surface conditions as 
\begin{align}
    \partial_t \phi_i &= - g \eta_i + p_D, \quad \text{on} \quad \Gamma^{\text{FS}} \times \mathcal{T}, \label{eq:DBC} \\
    \partial_t \eta_i &= \partial_z \phi_i + v_D, \quad \text{on} \quad \Gamma^{\text{FS}} \times \mathcal{T}, \label{eq:KBC}
\end{align}
with $\partial_t$ being the temporal gradient operator and $g$ is the constant gravitational acceleration. Furthermore, $v_D = v_D(x,y,t): \Gamma^{\text{FS}} \times \mathcal{T} \mapsto \R$ and $p_D = p_D(x,y,t): \Gamma^{\text{FS}} \times \mathcal{T} \mapsto \R$ are velocity and friction damping terms that are included to absorb outgoing waves following \cite{clamond2005efficient}, ultimately making the computational domain non-reflective. Those quantities are further specified in Section \ref{sec:absorption}.

\subsection{The radiation problem}\label{sec:radiation}

The radiation problems for $k = \{1,...,6\}$ represent the calm water responses to rigid-body motions in each of the $k$ directions. The pseudo-impulsive body boundary condition for the six radiation potentials, $\{\phi_k\}_{k=1}^6$, reads
\begin{equation}\label{eq:radiation_BC}
    \partial_n \phi_k = \partial_t x_k n_k, \quad \text{on} \quad \Gamma^{\text{body}},
\end{equation}
where $n_k = n_k(\boldsymbol{x}): \Gamma^{\text{body}} \mapsto \R$ is a generalized normal vector on $\Gamma^{\text{body}}$ defined by the unit normal vector, $\boldsymbol{n} = (n_x,n_y,n_z)$, and the position vector of a point on $\Gamma^{\text{body}}$, $\boldsymbol{r} = (r_x,r_y,r_z)$, by $n_k = \boldsymbol{n}$ for $k=\{1,2,3\}$ and $n_k = \boldsymbol{r} \times \boldsymbol{n}$ for $k = \{4,5,6\}$. Note that for additional generalized modes, equivalent normal vectors can be added. Moreover, $x_k = x_k(t): \mathcal{T} \mapsto \R$ is a pseudo-impulsive Gaussian displacement to be specified in Section \ref{sec:pseudo_impulse}. In the case of the classical impulsive formulation, $\partial_t x_k = \delta$ with $\delta = \delta(t)$ being the Dirac Delta function. A simple example of the body boundary condition can be seen in Figure \ref{fig:Pseudo_impulsive_boundary_conditions}.

As is apparent from Section \ref{sec:linear_potential_flow}, boolean values in terms of $\theta$ need to be specified for the radiation problem if symmetry boundaries are present. Those are listed in Table \ref{tab:Rad_BCs} for $k = \{1,...,6\}$ according to each of the solutions of $\phi_k$.

\begin{table}[h]
\footnotesize
 \caption{Boundary conditions for $\Gamma^{\text{sym,x}}$ and $\Gamma^{\text{sym,y}}$ in terms of $\theta$ for $k = \{1,...,6\}$ for $\phi_k$. $\theta = 0$ implies homogeneous Dirichlet type condition and $\theta = 1$ implies homogeneous Neumann type.}
 \vspace{2mm}
    \centering 
    \begin{tabular}{|c|c|c|c|c|c|c|} \hline
        $k$ & 1 & 2 & 3 & 4 & 5 & 6 \\ \hline
        $\Gamma^{\text{sym,x}} - (y=0)$ & $1$ & $0$ & $1$ & $0$ & $1$ & $0$\\ \hline
        $\Gamma^{\text{sym,y}} - (x=0)$ & $0$ & $1$ & $1$ & $1$ & $0$ & $0$ \\ \hline
    \end{tabular} \label{tab:Rad_BCs}
\end{table}

\subsubsection{Added mass and damping coefficients} 

For each of the radiation problems, $\phi_k$, the dynamic portion of the pressure can be computed from the linear Bernoulli equation as $p_k = \rho \partial_t \phi_k$, where $\rho$ is the constant water density, and $p_k = p_k(\boldsymbol{x},t): \Omega \times \mathcal{T} \mapsto \R$. Now, the radiation force in the $j$'th direction from a motion in the $k$'th direction, $F_{jk} = F_{jk}(t): \mathcal{T} \mapsto \R$, can be found via surface integration of the wetted body-boundary, $\Gamma^{\text{body}}$, as
\begin{equation}
    F_{jk} = \int_{\Gamma^{\text{body}}} p_k n_j d \Gamma, \quad \text{in} \quad \mathcal{T},
\end{equation}
with $n_j$ being the generalized normal vector as defined previously.

Let the added mass and damping coefficients be denoted by $a_{jk} = a_{jk}(\omega): \mathcal{H} \mapsto \R$ and $b_{jk} = b_{jk}(\omega): \mathcal{H} \mapsto \R$, respectively. Due to \cite{amini2017solving}, $a_{jk}$ and $b_{jk}$ can conveniently be computed as the Fourier transform of the radiation forces, $\F{F_{jk}}$, to the Fourier transform of the pseudo-impulsive displacement, $\F{x_{k}}$, as
\begin{equation}\label{eq:added_mass_and_damping}
    \omega^2 a_{jk} - i \omega b_{jk} = \frac{\F{F_{jk}}}{\F{x_{k}}}, \quad \text{in} \quad \mathcal{H}.
\end{equation}

\textbf{Infinite-frequency added mass}: At $\omega = \infty$, we define the infinite-frequency added mass as $a_{jk}^{\infty} = a_{jk}(\infty)$. This quantity can be computed by solving a single boundary value problem similar to the one presented above; however, with different boundary conditions on $\Gamma^{\text{FS}}$ and $\Gamma^{\text{body}}$ as
\begin{align}
\phi_k &= 0, \quad \text{on} \quad \Gamma^{\text{FS}}, \\
    \partial_n \phi_k &=  n_k, \quad \text{on} \quad \Gamma^{\text{body}}.
\end{align}

\subsection{The diffraction problem}\label{sec:diffraction}

The diffraction problem is the unknown scattered potential, $\phi_s$, due to fixed body interactions with a known incident potential, $\phi_0$. Hereby, the body boundary condition resembles an impermeability restriction which can be reformulated as
\begin{equation}\label{eq:diffraction_BC}
    \partial_n(\phi_0 + \phi_s) = 0 \quad \Longleftrightarrow \quad  \partial_n \phi_s = - \partial_n \phi_0, \quad \text{on} \quad \Gamma^{\text{body}}.
\end{equation}
The aim is to determine this condition in a pseudo-impulsive manner in the time domain as done for the radiation problem. In the frequency domain -- and in complex terms -- such a condition is straightforward to determine for a single wave frequency given the known incident velocity potential from linear wave theory \cite{newman1997marine}
\begin{equation}\label{eq:incident_wave_potential}
    \phi_0 = \operatorname{Re}\{\Psi e^{\i \omega t}\}, \quad \text{where} \quad \Psi = \frac{\i g}{\omega} \frac{\cosh(kh +kz)}{\cosh(kh)} e^{-\i k \alpha},
\end{equation}
with its associated free surface wave elevation
\begin{equation}\label{eq:incident_wave_elevation}
    \zeta_0 = \operatorname{Re}\{e^{-\i k \alpha} e^{\i \omega t}\}.
\end{equation}
Here the wave amplitude is assumed to be unity and $\zeta_0 = \zeta_0(x,y,t): \Gamma^{\text{FS}} \times \mathcal{T} \mapsto \R$. Please note how the incident wave is partitioned into a spatially varying part, $\Psi = \Psi(\boldsymbol{x}): \Omega \mapsto \R$, and a time phase, $e^{\i \omega t}$. The reasoning for this decomposition will become evident in Section \ref{sec:special}. Moreover, $\alpha = x \cos(\beta) + y \sin(\beta)$, with $\beta$ being the heading angle of the incoming waves defined relative to the positive $x$-axis in the $xy$-plane. 

From the classical linear theory of time-invariant dynamical systems, any velocity, $\nabla \phi_0$, can be determined by performing a convolution of the impulse response function, $\boldsymbol{K} = \boldsymbol{K}(\boldsymbol{x},t): \Omega \mapsto \R^3$, with any incident wave elevation, $\zeta_0(t) = \zeta_0(0,0,t)$, defined at the center of the horizontal plane, as
\begin{equation}
    \nabla \phi_0 (\boldsymbol{x},t) = \int_{-\infty}^{\infty} \boldsymbol{K}(\boldsymbol{x},\tau) \zeta_0(t-\tau) d \tau,
\end{equation}
where $\tau$ is an arbitrary integration variable. Now, inserting the know expressions for the incident wave potential and elevation from \eqref{eq:incident_wave_potential} and \eqref{eq:incident_wave_elevation} and rewriting
\begin{equation}
    \nabla \Psi(\boldsymbol{x}) e^{\i \omega t} =  \int_{-\infty}^{\infty} \boldsymbol{K}(\boldsymbol{x},\tau) e^{\i \omega (t-\tau)} d \tau = e^{\i \omega t} \int_{-\infty}^{\infty} \boldsymbol{K}(\boldsymbol{x},\tau) e^{-\i \omega \tau} d \tau,
\end{equation}
which simplifies by dividing through with $e^{\i \omega t}$ and using the exact definition of the radial-based Fourier transform from \eqref{eq:FT}
\begin{equation}\label{eq:some_eq}
    \nabla \Psi(\boldsymbol{x})  = \int_{-\infty}^{\infty} \boldsymbol{K}(\boldsymbol{x},\tau) e^{-\i \omega \tau} d \tau = \F{\boldsymbol{K}(\boldsymbol{x},\tau)}.
\end{equation}

Next -- as initially suggested in \cite{king1987time} -- an alternative approach would be to convolve the impulse response function, $\boldsymbol{K}(\boldsymbol{x},\tau)$, with another wave elevation, $\tilde{\zeta}_0$, that is not equal to \eqref{eq:incident_wave_elevation}, to determine the incident wave velocities, $\nabla \phi_0$. As for the radiation problem, we define this new $\tilde{\zeta}_0$ as a pseudo-impulsive Gaussian wave elevation with a tailored frequency content. See Section \ref{sec:pseudo_impulse} for the closed-form definition. Once again, we have the convolution integral
\begin{equation}
    \nabla \phi_0 (\boldsymbol{x},t) = \int_{-\infty}^{\infty} \boldsymbol{K}(\boldsymbol{x},\tau) \tilde{\zeta}_0(t-\tau) d \tau,
\end{equation}
which can be manipulated by taking Fourier transforms of the equation, which for convolution integrals equals multiplication in the frequency space as
\begin{equation}
    \F{\nabla \phi_0 (\boldsymbol{x},t)} = \F{\boldsymbol{K}(\boldsymbol{x},t)} \F{\tilde{\zeta}_0(t)}.
\end{equation}
This result is to be combined with \eqref{eq:some_eq} to give the final expression for the body boundary condition where the space and time dependencies have been dropped for simplicity
\begin{equation}\label{eq:diffraction_body_bc}
    \F{\nabla \phi_0} =  \nabla \Psi \F{\tilde{\zeta}_0} \quad \Longleftrightarrow \quad \nabla \phi_0 = \IF{\nabla \Psi \F{\tilde{\zeta}_0}}.
\end{equation}

All of the above are analytical expressions completely defined by linear potential flow wave theory and the pseudo-impulse to compute the velocity components necessary for \eqref{eq:diffraction_BC} combined with the unit body normal vector, $\boldsymbol{n} = (n_x,n_y,n_z)$. A simple example of the body boundary condition can be seen in Figure \ref{fig:Pseudo_impulsive_boundary_conditions} in Section \ref{sec:pseudo_impulse}.

\subsubsection{Exploiting geometric symmetry}\label{sec:special}

As previously noted, the geometrical features of the body and the fluid domain can allow for the utilization of symmetry boundaries, $\Gamma^{\text{sym},x}$ and $\Gamma^{\text{sym},y}$. For the radiation problems, the conditions on such boundaries are straightforward and intuitive; however, for the diffraction problem, one needs to pay special attention when incorporating symmetry conditions as the incident wave travels with a certain heading angle, $\beta$. Thus, we must make a spatial decomposition of the potential. First, considering the phase function of the incident wave in \eqref{eq:incident_wave_potential}, $e^{-\i k \alpha}$, which can be expanded using Euler's Formula as
\begin{equation}\label{eq:decomp}
    e^{-\i k \alpha} = \underbrace{\cos(k x c) \cos(k y s)}_{\text{SS}} - \i \underbrace{\cos(k x c) \sin(k y s)}_{\text{SA}} - \i \underbrace{\sin(k x c) \cos(k y s)}_{\text{AS}} - \underbrace{\sin(k x c) \sin(k y s)}_{\text{AA}},
\end{equation}
where the shorthand notation $c = \cos(\beta)$ and $s = \sin(\beta)$ have been used. This decomposition shows that the incident wave solution consists of four contributions, each characterized by whether they are symmetric or anti-symmetric around $(x=0)$ and $(y=0)$, respectively, hence the SS, SA, AS, and SS notation. Note how this result is solely related to the incident wave and not to the geometrical properties of $\Omega$.

Now, both the unknown scattered and known incident potentials are then decomposed as
\begin{equation}
    \phi_s = \phi_s^{\text{SS}} + \phi_s^{\text{SA}} + \phi_s^{\text{AS}} + \phi_s^{\text{AA}}, \quad \text{and} \quad \phi_0 = \phi_0^{\text{SS}} + \phi_0^{\text{SA}} + \phi_0^{\text{AS}} + \phi_0^{\text{AA}},
\end{equation}
where the former can be solved independently with suitable boundary conditions on the body boundary, $\Gamma^{\text{body}}$, and eventual symmetry boundaries $\Gamma^{\text{sym},x}$ and $\Gamma^{\text{sym},y}$. In this decomposed form, the latter conditions are close to trivial as listed in Table \ref{tab:Dif_BCs} in terms of the Boolean value, $\theta$, for each of the decomposed potentials.

The former condition -- on the body boundary, $\Gamma^{\text{body}}$ -- becomes of the following form
\begin{equation}\label{eq:body_bc_diffraction}
    \partial_n \phi_s^i = - \partial_n \phi_0^i = - \boldsymbol{n} \cdot \IF{\nabla \Psi^i \F{\tilde{\zeta}_0}}, \quad \text{on} \quad \Gamma^{\text{body}}, \quad \text{for} \quad i = \{\text{SS, SA, AS, AA}\}.
\end{equation}
We note that equation (2.30) \& (2.31) in \cite{amini2018pseudo} suggest a coupling between potentials in this boundary conditions, which is an unfortunate typographical error in that paper. A detailed deviation/explanation can be seen in Appendix \ref{sec:appendix_1}.

\textbf{Remarks on having one symmetry boundary}: We assume the possibility of two symmetry boundaries in the above. In the case of having only one of those -- either $(x=0)$ or $(y=0)$ -- \eqref{eq:decomp} reduces to having only a symmetrical and an anti-symmetrical term. Ultimately yielding $\phi_s = \phi_s^{\text{S}} + \phi_s^{\text{A}}$ and $\phi_0 = \phi_0^{\text{S}} + \phi_0^{\text{A}}$.

\begin{table}[h]
\footnotesize
 \caption{Boundary conditions for $\Gamma^{\text{sym,x}}$ and $\Gamma^{\text{sym,y}}$ in terms of $\theta$ for $\left (\phi_s^{\text{SS}},\phi_s^{\text{SA}}, \phi_s^{\text{AS}}, \phi_s^{\text{AA}} \right)$. $\theta = 0$ implies homogeneous Dirichlet type condition and $\theta = 1$ implies homogeneous Neumann type.}
 \vspace{2mm}
    \centering 
    \begin{tabular}{|c|c|c|c|c|} \hline
         & SS & SA & AS & AA \\ \hline
        $\Gamma^{\text{sym,x}} - (y=0)$ & $1$ & $0$ & $1$ & $0$ \\ \hline
        $\Gamma^{\text{sym,y}} - (x=0)$ & $1$ & $1$ & $0$ & $0$ \\ \hline
    \end{tabular} \label{tab:Dif_BCs}
\end{table}

\subsubsection{Wave excitation forces}

The scattered wave forces in the $j$'th direction, $F_{s,j}$, is computed via integration of the scattered wave pressure, $p_s = p_s(\boldsymbol{x},t): \Omega \times \mathcal{T} \mapsto \R$, over the wetted body surface as
\begin{equation}\label{eq:scattered_force}
    F_{s,j} =  \int_{\Gamma^{\text{body}}} p_s n_j d \Gamma, \quad \text{in} \quad \mathcal{T},
\end{equation}
where $n_j$ is the generalized normal vector as for the radiation forces. Following \cite{amini2018pseudo}, the scattered wave excitation force can be computed by
\begin{equation}\label{eq:wave_exitation}
    X_{s,j} = \frac{\F{F_{s,j}}}{\F{\tilde{\zeta}_0}}, \quad \text{in} \quad \mathcal{H}.
\end{equation}

The complete diffraction wave excitation force also consists of an incident wave contribution known as the Froude-Krylov force, \cite{faltinsen1993sea}, such that $X_{D,j} = X_{0,j} + X_{s,j}$. In the frequency domain, the incident wave pressure, $p_0 = p_0(\boldsymbol{x},t): \Omega \times \mathcal{T} \mapsto \R$, can be integrated over $\Gamma^{\text{body}}$ as in \eqref{eq:scattered_force} using
\begin{equation}
    p_0 = \rho g \frac{\cosh(kh +kz)}{\cosh(kh)} e^{-\i k \alpha}.
\end{equation}

\textbf{Remark on the symmetry decomposition}: Assuming two symmetry boundaries. Due to the decomposition, the different excitation forces solely depend on the solution that matches the symmetry or anti-symmetry of the problem. For completeness,
\begin{itemize}
    \item SS forces: Heave from $X_{D,3} = X_{D,3}(p_0^{\text{SS}},p_s^{\text{SS}})$.
    \item SA forces: Sway from $X_{D,2} = X_{D,2}(p_0^{\text{SA}},p_s^{\text{SA}})$ and roll from $X_{D,4} = X_{D,4}(p_0^{\text{SA}},p_s^{\text{SA}})$.
    \item AS forces: Surge from $X_{D,1} = X_{D,1}(p_0^{\text{AS}},p_s^{\text{AS}})$ and pitch from $X_{D,5} = X_{D,5}(p_0^{\text{AS}},p_s^{\text{AS}})$.
    \item AA forces: Yaw from $X_{D,6} = X_{D,6}(p_0^{\text{AA}},p_s^{\text{AA}})$.
\end{itemize}
Similarly, in the case of one symmetrical boundary only having pure symmetrical or anti-symmetrical contributions.

\subsection{The pseudo-impulsive Gaussian displacement and wave elevation}\label{sec:pseudo_impulse}

As is evident from Section \ref{sec:radiation} and Section \ref{sec:diffraction}, the system dynamics -- or the frequency range hereof -- depends on the pseudo-impulse ($x_k$ for the radiation problems and $\tilde{\zeta}_0$ for the diffraction problem). Following \cite{amini2017solving,amini2018pseudo,visbech2023spectral} and the initial work by \cite{king1987time}, we define a Gaussian pseudo-impulse function, $g = g(t): \mathcal{T} \mapsto \R$, of unit height by
\begin{equation}\label{eq:Gaussian}
 g(t) = e^{-2 \pi^2 s^2 (t-t_0)^2}.
\end{equation}
With this, the width -- or variance -- of the function is controlled by $s$, and the peak location is determined as $t_0 = \sqrt{\frac{\log(\epsilon)}{-2 \pi^2 s^2}}$, such that the function is practically zero at $t=0$ when $\epsilon$ is chosen to be sufficiently small. Smaller values of $s$ increase the frequency span of the final solution and vice versa. The discrete considerations of this function -- and how $s$ is defined as a function of the spatial resolution -- are omitted in this paper; however, a detailed summary is outlined in \cite{visbech2023spectral}. We use this impulse for all the radiation problems as $x_k(t) = g(t)$ and for the diffraction problem as $\tilde{\zeta}_0(t) = g(t)$. Note that analytical expression can be obtained for $\hat{g}(\omega) = \F{g(t)}$ to be used efficiently in the computation of added mass and damping coefficients, the scattered body boundary condition, and wave excitation forces.

\begin{figure}[h]
    \centering
    \includegraphics[width=0.70\textwidth]{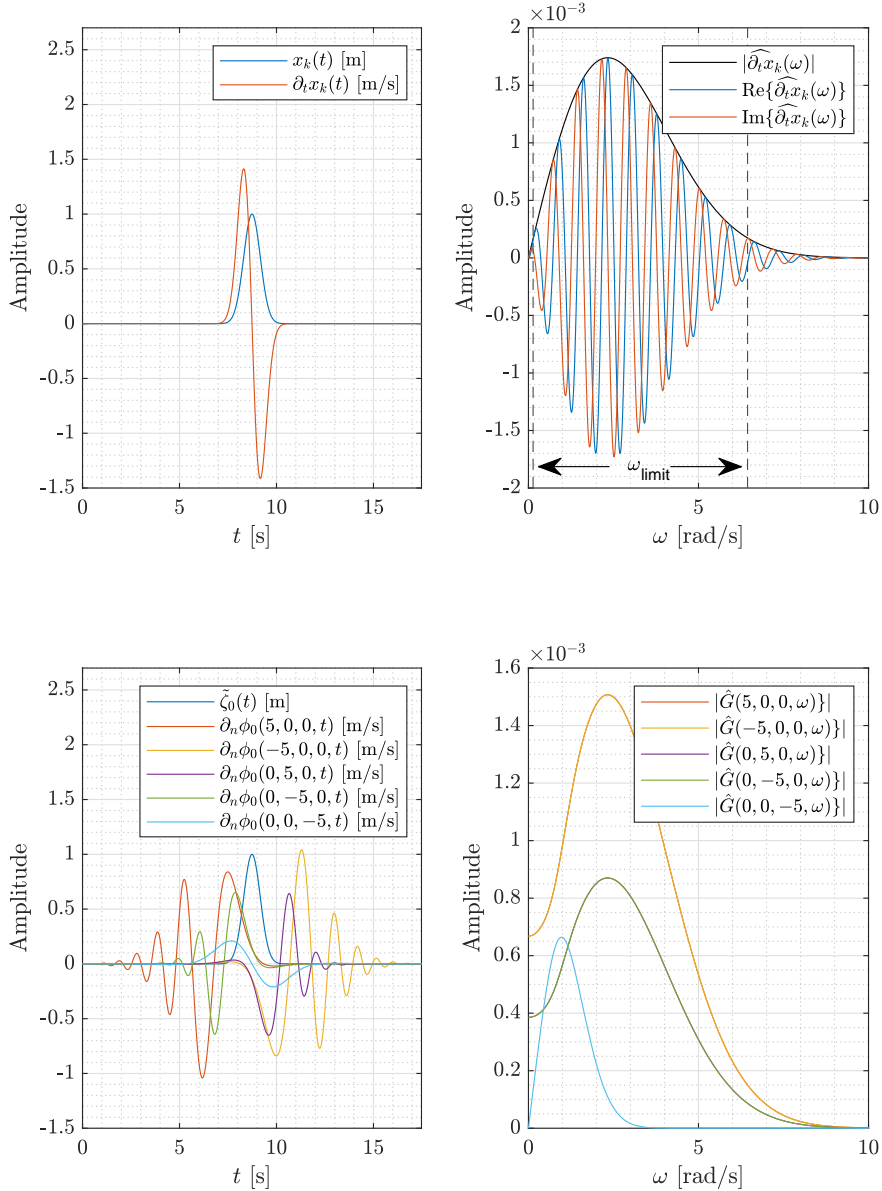}
    \caption{A simple example of the pseudo-impulsive radiation and diffraction body boundary condition in the time and frequency domain of a sphere with radius $R = 5 $ [m] with $\beta = \frac{7}{8} \pi$.
    \textit{Top-left}: The pseudo-impulsive displacement, $x_k(t)$, and its derivative, $\partial_t x_k(t)$ for the radiation problem in the time domain.
    \textit{Top-right}: The pseudo-impulsive velocity in the frequency domain, $\widehat{\partial_t x}_k(\omega)$, and various measures hereof. Including $\omega_{\text{limit}}$, which contains all frequencies that have more than 10\% of the maximum absolute energy.
   \textit{Bottom-left}: The pseudo-impulsive wave elevation, $\tilde{\zeta}_0(t)$, and the associated diffraction body boundary condition, $\partial_n \phi_0(x,y,z,t)$, on different locations of the sphere.
    \textit{Bottom-right}: The diffraction body boundary condition in the frequency domain as $\hat{G}(x,y,z,\omega) = \F{\partial_n \phi_0(x,y,z,t)}$.}
    \label{fig:Pseudo_impulsive_boundary_conditions}
\end{figure}

In Figure \ref{fig:Pseudo_impulsive_boundary_conditions}, a simple example is given of the pseudo-impulsive radiation and diffraction body boundary condition for a sphere of radius $R = 5$ [m] with wave approaching with an angle $\beta = \frac{7}{8} \pi$. Starting from the \textit{top-left} plot, the pseudo-impulsive displacement -- given by \eqref{eq:Gaussian} -- and its time derivative (a pseudo-impulsive velocity) is shown in the time domain. On the \textit{top-right} plot, the associated frequency domain of the velocity is pictured. Note how the frequency span of the impulse can be widened by decreasing the width of the displacement in the time domain. Also, in the frequency domain, we define a span denoted by $\omega_{\text{limit}}$ that constitutes frequencies with more than 10\% of the maximum absolute energy. This span will be the range where one can expect a well-behaved solution. For the numerical results in Section \ref{sec:numerical_results}, all frequency-based results will be plotted within $\omega_{\text{limit}}$. Moving on to the \textit{bottom-left} figure, the pseudo-impulsive wave elevation and associated diffraction boundary condition at various locations of the sphere (at the front, the back, on the sides, and on the bottom) are plotted. The associated Fourier transform of those conditions is shown on the last and \textit{bottom-right} plot. Note how the condition at the front and back impose the same frequency content as the ones on the sides.

\subsection{Wave absorption}\label{sec:absorption}

A significant and very important task when building numerical wave tanks is obtaining a computational domain free of reflections, as such reflections would influence the final frequency content of the solution. Several possible approaches and techniques ensure a domain with those properties. For this work, we combined: i) grid stretching and ii) pressure and velocity damping zones.

\textbf{Grid stretching:} By stretching the initial mesh, the computational domain can be expanded to sufficient size without increasing the computational cost. Emphasis is put on \textit{gentle}, as too rapid changes in grid size have been seen to induce reflections if dynamics are present on the pre-stretched mesh. If the mesh subject to grid stretching primarily contains the zero solution, one can employ strong stretching without introducing any reflections. This trick can be used to capture the low-frequency content of the radiation problem.

\textbf{Pressure and velocity damping zones:} As is apparent from the free surface boundary conditions in \eqref{eq:DBC} and \eqref{eq:KBC}, two additional terms -- $p_D$ and $v_D$ -- have been added to absorb outgoing waves. Those pressure and friction absorption terms are incorporated as damping zones following \cite{clamond2005efficient} using spatially varying functions $c_p = c_p(x,y): \Gamma^{\text{FS}} \mapsto \R$ and $c_v = c_v(x,y): \Gamma^{\text{FS}} \mapsto \R$. Moreover, $c_p$ and $c_v$ are, per definition, zero outside the damping zone and resemble great similarities with the function used in e.g., classical relaxations zones by \cite{larsen1983open}. For this work, we employ a Gaussian-type function that smoothly transitions and has a maximum of $2\pi$ in the middle of the damping zone.

The velocity damping terms are defined as
\begin{equation}
    v_D = -c_v \eta_i, \quad \text{on} \quad \Gamma^{\text{FS}},
\end{equation}
and for the pressure damping term, we have the following scalar Poisson equation
\begin{equation}\label{eq:pressure_damping_term}
    \nabla^2 p_D =  \nabla \cdot \left ( c_p \nabla \phi_i \right ), \quad \text{on} \quad \Gamma^{\text{FS}},
\end{equation}
with $p_D=0$ and $\partial_n p_D = 0$ conditions according to the design of the damping zone. Please note that in relation to wave absorption, $\nabla = (\partial_x,\partial_y)$ is the Cartesian differential operator on $\Gamma^{\text{FS}}$ (the $xy$-plane) and $\partial_n$ is only the normal direction derivative on this plane. In the following, $\nabla$ is reconsidered as the classical three-dimensional operator in $\Omega$ until further notice. Moreover, the pressure term differs slightly from the one presented in \cite{clamond2005efficient}, as the gradient has been taken of the original equation followed by the divergence. Hereby, we obtain the scalar Poisson equation on $\Gamma^{\text{FS}}$. The discrete solution approach to this problem is discussed in Section \ref{sec:computing_pressure_damping_term}.
\section{Numerical Methods}\label{sec:discretization_methods}

The classical \textit{method of lines} (MOL) approach is adopted to discretize the problem in space via the spectral element method while keeping it continuous in time. This spatial discretization of $\Omega$ will be further outlined in Section \ref{sec:spectral_element_method}. In Section \ref{sec:hybrid_meshes}, we highlight meshes and the generation of so-called \textit{hybrid} meshes, which provides a stable geometrical basis for the spectral element-based free surface model. Moreover, the remaining discretization tasks are outlined briefly in the following:

\begin{itemize}
    \item \textit{Time integration}: Due to the MOL, we transform \eqref{eq:DBC} and \eqref{eq:KBC} into a semi-discrete system of coupled ordinary differential equations. The discrete-time integration hereof is performed by a traditional explicit four-stage fourth-order Runge-Kutta method (ERK4) with a suitable uniform time-step governed by a CFL condition: $\Delta t = C \Delta x_{\min} u_{\max}^{-1}$, where $u_{\max} = \sqrt{g h}$ is chosen to be the asymptotic shallow water limit, $\Delta x_{\min}$ is a global measure of the smallest node-spacing on $\Gamma^{\text{FS}}$, and $0 < C \leq 1$. This time step is intentionally designed conservatively such that high-frequency waves are well-discretized in time.

    \item \textit{Time differentiation}: As outlined and detailed in \cite{visbech2023spectral}, time derivatives, $\partial_t$, will be approximated to fourth order in accuracy -- to accommodate the ERK4 -- using finite differences stencils. This approach has been used to compute the dynamic pressure, $p_i$, from the scattered and radiation potential and the pseudo-impulsive velocities for the radiation problems in Section \ref{sec:radiation}. 

    \item \textit{Fourier and inverse Fourier transformations}: The various Fourier transforms required to compute the added mass and damping coefficients in \eqref{eq:added_mass_and_damping} and the scattered wave excitation force in \eqref{eq:wave_exitation}, will be evaluated discretely via fast Fourier transforms. The inverse operation is performed similarly for the inverse operator in connection with the evaluation of the body boundary condition for the diffraction problem in \eqref{eq:diffraction_body_bc}.
    
\end{itemize}

\subsection{A spectral element approximation}\label{sec:spectral_element_method}

Considering the spatial domain, $\Omega$, we apply a complete element-tessellation of the volume and its surfaces. This partitioning is performed using so-called \textit{hybrid} meshes (see Section \ref{sec:hybrid_meshes}), where a series of prism and tetrahedral elements are arranged in a boundary-conform non-overlapping fashion. The complete number of elements, $N_{\text{elm}}$, is taken as the sum of the number of the two individual element types and therefore defined as $N_{\text{elm}} = N_{\text{elm}}^{\text{pri}} + N_{\text{elm}}^{\text{tet}}$. With this, $\Omega \simeq \bigcup_{n = 1}^{N_{\text{elm}}} \mathcal{E}^n$, where $\mathcal{E}^n$ is the $n$'th element. For the SEM, we define a space of continuous and piece-wise polynomial functions with degree at most $P$
\begin{equation}
    \mathcal{V}^P = \{v \in C^0(\Omega); \forall n \in \{1,..., N^{\text{elm}} \}, v \mid_{\mathcal{E}^n} \in {\rm I\!P}^{P}\}.
\end{equation}

Now, using the geometrical representation of $\Omega$, a global representation of any velocity potential, $\phi$, reads
\begin{equation}\label{eq:geo_approx}
    \phi \simeq \bigoplus_{n = 1}^{N_{\text{elm}}} \phi^n, \quad \text{where} \quad \phi^n = \phi(\overline{\boldsymbol{x}}) \quad \text{for} \quad \overline{\boldsymbol{x}} \in \mathcal{E}^n,
\end{equation}
where $\phi$ is assumed to be geometrical time-invariant (due to linearity) such that no re-meshing strategy is needed. Moreover, $\phi^n$ denotes the local solution of $\phi$ on $\mathcal{E}^n$. The element-defined coordinate, $\overline{\boldsymbol{x}}$, is given as $N_{\text{vp}}$ Gauss-Lobatto-type distributed volume points (quadrature points), such that $\overline{\boldsymbol{x}} = \{x_i,y_i,z_i \}_{i=1}^{N_{\text{vp}}}$. Now, for each local solution, we apply a modal-nodal polynomial representation, e.g., see \cite{karniadakis2005spectral}, as
\begin{equation}\label{eq:sol_approx}
         \phi^n \approx \sum_{i=1}^{N_{\text{vp}}} \hat{\phi}^n_{i} \psi_{i}(\overline{\boldsymbol{x}}) = \sum_{j=1}^{N_{\text{vp}}} \phi_{j}^n h_{j}(\overline{\boldsymbol{x}}), \quad \text{for} \quad \overline{\boldsymbol{x}} \in \mathcal{E}^n.
\end{equation}
Here the former sum constitutes the modal representation, and the latter is the nodal. For the modal representation, $\{\hat{\phi}^n_{i} \}_{i=1}^{N_{\text{vp}}}$ is a set of expansion coefficients connected with the three-dimensional modal Jacobi-type basis functions, $\psi_{i}(\overline{\boldsymbol{x}}) \in \mathcal{V}^P$, up to order $P$. For the nodal representation, $\{\phi_{j}^n \}_{j=1}^{N_{\text{vp}}} = \phi^n(\{x_j,y_j,z_j \}_{j=1}^{N_{\text{vp}}})$ is a set of nodal solution values at the quadrature points, which connects to the three-dimensional nodal Lagrange polynomial basis function, $h_{j}\in \mathcal{V}^P$, up to order $P$.

Following \cite{hesthaven2007nodal}, an iso-parametric solution strategy is adopted, where we define a Cartesian reference domain, $\mathcal{R}$, spanned by $\boldsymbol{r} = (r,s,t)$. Each element can then be mapped back and forth through the mapping $\overline{\Psi} = \overline{\Psi}(\boldsymbol{r}): \mathcal{E}^n \mapsto \mathcal{R}$ and its inverse such that $\overline{\boldsymbol{x}} = \overline{\Psi}(\boldsymbol{r})$ and $\boldsymbol{r} = \overline{\Psi}^{-1}(\overline{\boldsymbol{x}})$. An orthonormal basis is constructed on $\mathcal{R}$, leading to well-behaved spatial operators for differentiation, integration, and interpolation. For a more detailed run-down of the SEM approximation, consider the key literature in \cite{karniadakis2005spectral,hesthaven2007nodal} or in the more applied context \cite{engsigkarup2016stabilised,xu2018spectral,visbech2023spectral}.

\subsubsection{Weak Galerkin spectral element discretization}

The SEM is a Galerkin-based method, which implies -- among other things -- that the governing equations are satisfied in a weak sense. The weak formulation of the pseudo-impulsive problem is obtained by: i) multiplying the Laplace equation in $\Omega$ with a test function, $v = v(\boldsymbol{x}): \Omega \mapsto \R$, ii) integrating over the domain, and iii) apply Green's 1st identity. Doing so, one gains less regularity on the solution as $\phi \in C^1(\Omega)$, but also naturally incorporates Neumann-type boundary conditions. The final weak form of the problem yields
\begin{equation}\label{eq:weak}
    \int_{\Omega} \nabla \phi \cdot \nabla v ~d \Omega = \int_{\Gamma^{\text{body}}} q v ~d \Gamma,
\end{equation}
where $\phi = \phi_k$ and $q = \partial_t x_k n_k$ for the radiation problems and $\phi = \phi_s$ and $q = - \partial_n \phi_0$ for the diffraction problem.

Choosing the test functions, $v$, to equal the nodal basis functions, $h$, a nodal Galerkin discretization is adopted. Next, by inserting the solution approximation from \eqref{eq:geo_approx} and \eqref{eq:sol_approx} with the geometrical approximation into the weak formulation in \eqref{eq:weak}, a discrete linear system of equations is constructed as
\begin{equation}\label{eq:linear_system_of_equations}
    \mathcal{A} \boldsymbol{\phi} = \boldsymbol{b}, \quad \text{where} \quad \mathcal{A} \in \R^{N_\text{DOF} \times N_\text{DOF} }, \quad \text{and} \quad \{ \boldsymbol{\phi},\boldsymbol{b} \} \in \R^{N_\text{DOF}}.
\end{equation}
Here, $N_{\text{DOF}}$ denotes the number of degrees of freedom (DOF) in the discrete scheme. Moreover, the system matrix is represented as $\mathcal{A}$, the system vector as $\boldsymbol{b}$, and the global solution vector as $\boldsymbol{\phi}$. Local three-dimensional element contributions are utilized to construct $\mathcal{A}$, resulting in a sparse and symmetric matrix. Similarly, $\boldsymbol{b}$ is formed by incorporating local two-dimensional element contributions. To impose the inhomogeneous Dirichlet boundary data, where the test functions vanish on the boundaries of the Dirichlet surface, i.e., $v(\boldsymbol{x}) = 0$ for $\boldsymbol{x} \in \Gamma^{\text{FS}}$, modifications are made to the system matrix and vector during a post-assembling step.

\subsubsection{Computing the pressure damping term}\label{sec:computing_pressure_damping_term}

To obtain the pressure damping term in \eqref{eq:pressure_damping_term} from Section \ref{sec:absorption}, one needs to solve a two-dimensional scalar Poisson problem. The weak formulation of this problem reads
\begin{equation}
      \int_{\Gamma^{\text{FS}}} \nabla p_D \cdot \nabla v ~d \Gamma = - \int_{\Gamma^{\text{FS}}} \nabla \cdot \left ( c_p \nabla \phi_i \right ) v ~d \Gamma,
\end{equation}
which gives rise to another linear system of equations, $\mathcal{A} \boldsymbol{p_D} = \boldsymbol{b}$, where homogeneous Dirichlet boundary data are incorporated as a post-assembling step as well. Recall that $\nabla$ in the case of wave absorption is the two-dimensional gradient operator in the $xy$-plane.

\begin{figure}[t] 
    \centering
     \fbox{\includegraphics[scale=0.5]{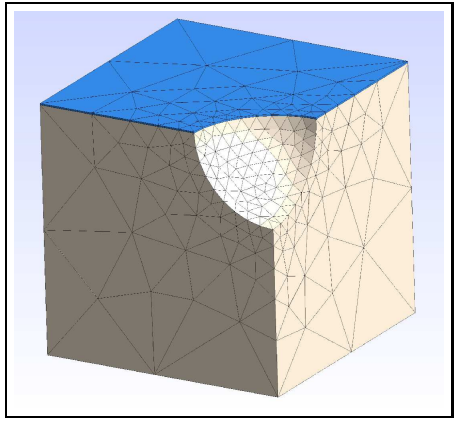}}
    \caption{Hybrid mesh of a floating sphere with two symmetry boundaries. Snapshot from Gmsh.}
    \label{fig:hybrid_mesh}
\end{figure}

\subsection{Hybrid meshes}\label{sec:hybrid_meshes}

Numerical models that employ discrete geometrical representation -- such as finite and spectral element methods -- rely heavily on the quality of the mesh. Features like mesh refinement (the possibility to refine -- or grade -- the mesh towards specific parts of the domain), unstructuredness of the mesh (having a dynamic mesh with non-structured elements), and mesh generation are of massive importance. Another feature that relates to solving free surface flows using the SEM is mesh-induced temporal instabilities. In the late 1990s, \cite{robertson1999free} pointed out stability issues when solving a nonlinear free surface potential flow problem using structured triangular elements. The primary reason for this was the asymmetrical spatial discretization of the free surface. Two decades later, \cite{engsigkarup2016stabilised} succeeded in stabilizing the problem by employing structured quadrilateral elements. Then, in \cite{engsigkarup2019mixed}, a composite -- or \textit{hybrid} -- approach was used with a single layer of structured quadrilateral elements at the free surface (for symmetry) followed by unstructured triangular elements underneath. At last, in the three-dimensional setting, \cite{engsigkarup2018spectral} employed structured (in the $z$-direction) prism elements.

\subsubsection{The ultimate meshing approach for incompressible and inviscid free surface flows}

In terms of geometrical flexibility, the three-dimensional simplex (the tetrahedral) is arguably one of the most versatile types of elements. In the case of SEM and PF problems with the presence of floating bodies, such meshes are difficult to construct with the symmetry restriction and will -- with great possibility -- be highly unstable. Therefore, we propose the ultimate meshing approach for incompressible and inviscid free surface flow for SEM by exploiting -- what we call \textit{hybrid} meshes -- as
\begin{itemize}
    \item A single and arbitrarily thin -- yet of finite thickness -- top layer of prism elements. To cope with the symmetry requirement, those elements are arranged unstructured in the $xy$-plane and structured in the $z$-direction.
    \item Followed from underneath is a completely unstructured array of tetrahedral elements to capture any complex geometrical features of the body and other boundaries. This set of elements will constitute the majority of the complete discrete representation and allow for any key meshing feature.
\end{itemize}

Hybrid meshes for this paper are constructed using Gmsh, \cite{geuzaine2009gmsh}. Figure \ref{fig:hybrid_mesh} shows a meshing example of a floating sphere. Note how the mesh is graded towards the sphere and less towards the remaining boundaries. 

\section{Validation}\label{sec:validation_and_computational_properties}

 \begin{figure}[b]
    \begin{minipage}[t]{.45\textwidth}
        \centering
        \includegraphics[scale=1]{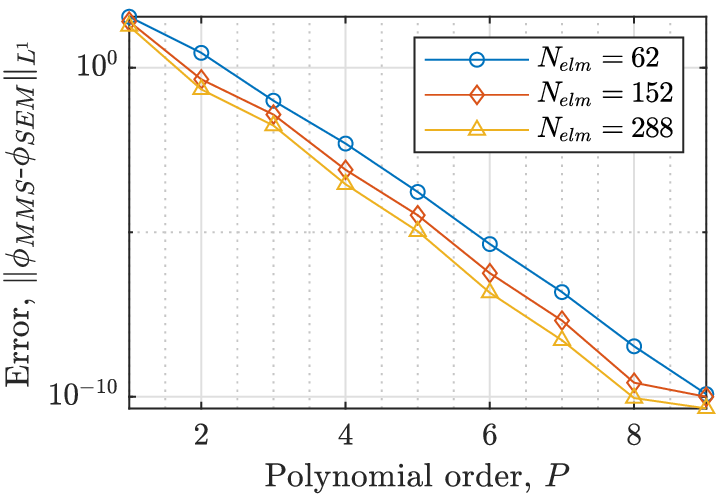}
        \subcaption{Polynomial order, $P = \{1,...,9\}$, on three different meshes,  $N_{\text{elm}} = \{62, 152, 288\}$.}
    \end{minipage}
    \hfill
    \begin{minipage}[t]{.45\textwidth}
        \centering
        \includegraphics[scale=1]{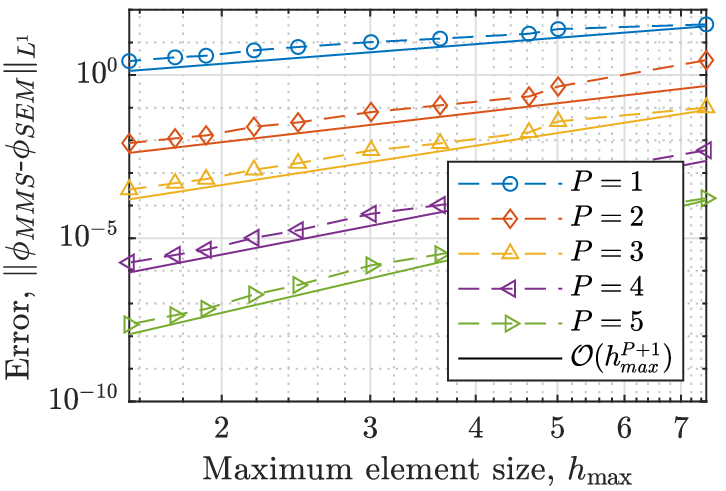}
        \subcaption{Maximum element size, $h_{\max}$, for meshes with $N_{\text{elm}} \in [62 ~;~ 4,872]$ and polynomial order in $P = \{1,...,5\}$.}
    \end{minipage}  
    \caption{Convergence studies.}
    \label{fig:Convergence}
\end{figure}

In this section, we validate the proposed numerical solver by carrying out spatial convergence (algebraic and spectral) studies. First, consider a computational domain of a floating sphere with radius, $R = 5$, that is placed at $\boldsymbol{x} = (0,0,0)$. Two symmetry boundaries are employed; thus, only a quarter of the sphere is represented, and all three domain axes are of length $1.5R$. Meshing is of the hybrid type with elements uniformly in size. To verify the discrete Laplace operator, the \textit{method of manufactured solutions} (MMS) is used, e.g., see \cite{roache2019method}. Hereby, the Laplace problem is recast into a Poisson problem with a known infinitely smooth MMS-solution, $\phi_{\text{MMS}} \in C^{\infty}(\Omega)$, of trigonometric functions. Dirichlet boundary conditions are imposed on the top boundary, and Neumann conditions are enforced on the remaining boundaries. Ultimately, we want to show and confirm $p$- and $h$-convergence. That is exponential/spectral decay in error for a fixed mesh under the increase of polynomial order, $P$, and algebraic (of fixed polynomial order) under the decrease in element size, $h$, respectively. The error norm for the preceding convergence studies is chosen to be the exact global $L_1$-norm, $\| \phi_{\text{MMS}}-\phi_{\text{SEM}}\|_{L_1} = \int_{\Omega} |\phi_{\text{MMS}}(\boldsymbol{x})-\phi_{\text{SEM}}(\boldsymbol{x}) |^1 ~d \Omega$ by the use of a global integration matrix.

\subsection{$\boldsymbol{p}$- and $\boldsymbol{h}$-convergence}

For the $p$-convergence study, we consider three different meshes given by the total number of elements, $N_{\text{elm}} = \{62, 152, 288\}$. For each mesh, we apply polynomial order in the range $P = \{1,...,9\}$, as this is the mesh generation limit for Gmsh. The combination of meshes and polynomial orders gives systems with the sizes $N_{\text{DOF}} \in [34 ~;~ 44,491]$. From Figure \ref{fig:Convergence} (a), spectral convergence can be confirmed visually in the semi-logarithmic plot. Errors reach below $10^{-10}$ and drop about 12 orders of magnitude over the nine polynomial orders. At this point ($P = 9$), we expect round-off errors to have been accumulated and hereby begin to influence the solution. The finest meshes yield smaller errors than the coarser ones, which indicates $h$-convergence. For that study, we consider polynomial order in the range $P = \{1,...,5\}$, combined with meshes ranging within $N_{\text{elm}} \in [62 ~;~ 4,872]$, yielding system sizes in the range $N_{\text{DOF}} \in [34 ~;~ 118,701]$. From Figure \ref{fig:Convergence} (b), algebraic convergence can be observed with order $\O{h^{P+1}}$.

\section{Numerical Results}\label{sec:numerical_results}

The upcoming section presents the spectral element solution to the three-dimensional linear pseudo-impulsive radiation and diffraction problem. First, we show \textit{proof-of-concept} cases in terms of a floating sphere and a floating box to further validate the legitimacy of the numerical model. Secondly, a more complex case is considered in terms of an oscillating water column (OWC). This latter case considers the use of special boundaries and generalized modes. 

In the following, results for added mass and damping coefficient, $a_{jk}$ and $b_{jk}$, and wave excitation forces, $X_{0,j}$ and $X_{s,j}$, will be shown and compared with associated WAMIT results. These physical quantities are non-dimensionalized according to \cite{lee2006wamit} as
\begin{equation}
    \bar{a}_{jk} = \frac{a_{jk}}{\rho L^n}, \quad  \bar{b}_{jk} = \frac{b_{jk}}{\rho L^n \omega}, \quad \text{and} \quad \bar{X}_j = \frac{X_j}{\rho g L^m},
\end{equation}
where $n = 3$ for $(j,k) = \{1,2,3\}$, $n = 4$ for $j = \{1,2,3\}$ and $k = \{4,5,6\}$, $n = 4$ for $j = \{4,5,6\}$ and $k = \{1,2,3\}$, and $n = 5$ for $(j,k) = \{4,5,6\}$. Moreover, $m = 2$ for $j = \{1,2,3\}$ and $m = 3$ for $j = \{4,5,6\}$. Here, $L$ is the length scale of the structure in question. Also, a non-dimensional radian frequency is defined by $\overline{\omega} = \omega \sqrt{L/g}$. The frequency-based quantities will be plotted within the $\omega_{\text{limit}}$ as detailed and highlighted in relation to Figure \ref{fig:Pseudo_impulsive_boundary_conditions}.

\subsection{A floating sphere}\label{sec:sphere}

For the first numerical benchmark test, we consider a floating sphere of radius, $R = 5$ [m], such that the characteristic length scale corresponds to $L = R$. The water depth of the simulation is set to $h = 5 R$, and the fluid domain exploits two symmetry boundaries for efficient computation with a hybrid mesh that is constructed in a graded fashion towards the body boundary, $\Gamma^{\text{body}}$. The polynomial order is chosen to be $P=3$, and the quarter of the sphere is meshed using $174$ surface elements. Moreover, the total number of volumetric elements are $N_{\text{elm}} = 1,629$, yielding a system with $N_{\text{DOF}} = 11,787$.

\subsubsection{Added mass and damping coefficients} 

\begin{figure}[t]
    \begin{minipage}[t]{.45\textwidth}
        \centering
        \includegraphics[width=0.85\textwidth]{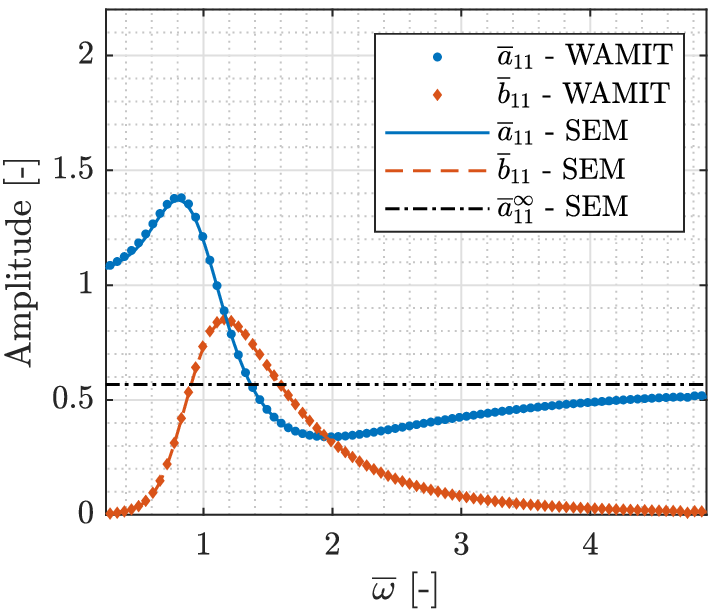}
        \subcaption{Surge-Surge}
    \end{minipage}
    \hfill
    \begin{minipage}[t]{.45\textwidth}
        \centering
        \includegraphics[width=0.85\textwidth]{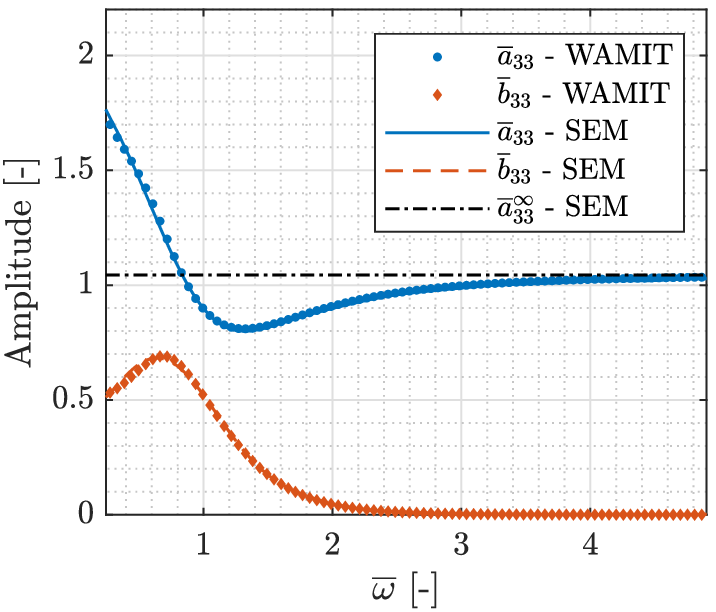}
        \subcaption{Heave-Heave}
    \end{minipage}  
    \caption{Added mass and damping coefficients, $\bar{a}_{jk}$ and $\bar{b}_{jk}$, for a floating sphere of radius $R$ on finite depth, $h = 5R$. Comparison between the proposed SEM model and WAMIT results.}
    \label{fig:sphere_a_jk_b_jk}
\end{figure}

The radiation problem is solved in surge and heave. With this, we consider the added mass and damping coefficient in relation hereof, namely, $\bar{a}_{11}$, $\bar{a}_{33}$, $\bar{b}_{11}$, and $\bar{b}_{33}$. The coefficients are shown in Figure \ref{fig:sphere_a_jk_b_jk}. From the figure, a good visual agreement can be seen when compared with the simulation results obtained by WAMIT. Also, the added mass clearly convergence towards the infinite added mass. 

\subsubsection{Wave excitation forces}

\begin{figure}[t]
    \begin{minipage}[t]{.45\textwidth}
        \centering
        \includegraphics[width=0.85\textwidth]{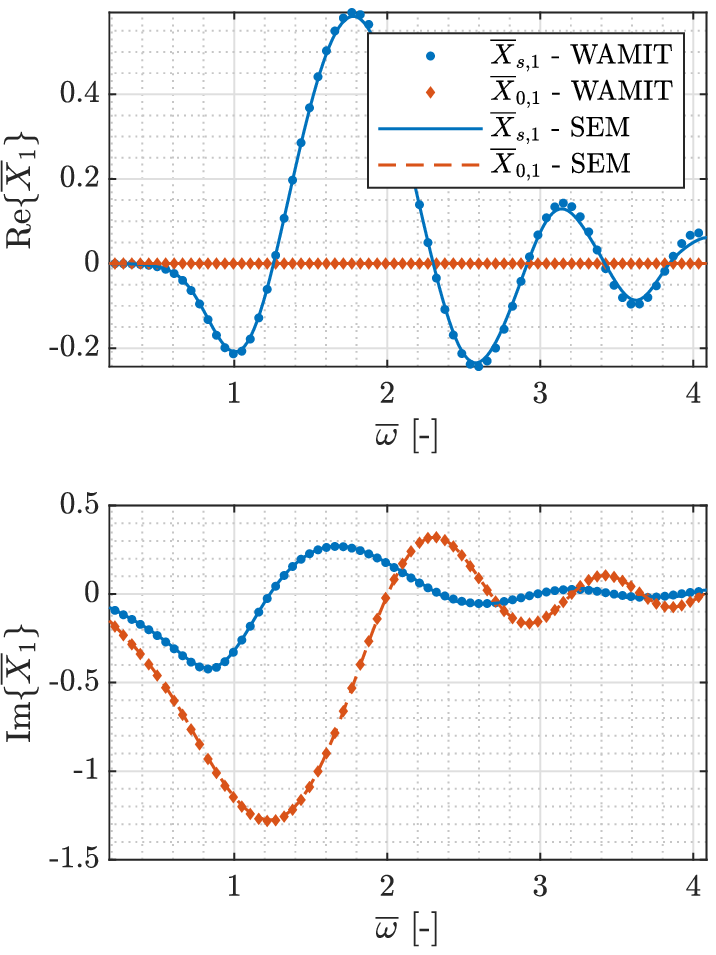}
        \subcaption{Surge}
    \end{minipage}
    \hfill
    \begin{minipage}[t]{.45\textwidth}
        \centering
        \includegraphics[width=0.85\textwidth]{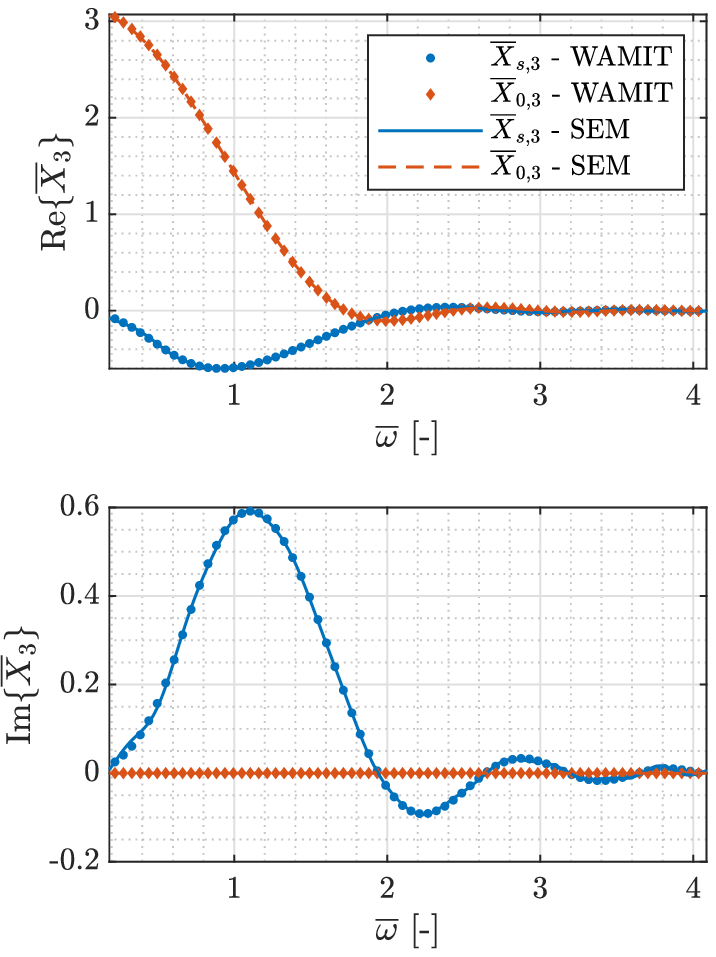}
        \subcaption{Heave}
    \end{minipage}  
    \caption{Wave excitation forces, $\bar{X}_{0,j}$ and $\bar{X}_{s,j}$ for a floating sphere of radius $R$ on finite depth, $h = 5R$, with heading angle $\beta = 150$ degrees. Comparison between the proposed SEM model and WAMIT results.}
    \label{fig:sphere_X_j}
\end{figure}

Next, the diffraction problem is solved, where we consider the wave excitation forces for the scattered and the incident wave (the Froude-Krylov forces). This is done for surge, $\bar{X}_{0,1}$ and $\bar{X}_{s,1}$, and heave, $\bar{X}_{0,3}$ and $\bar{X}_{s,3}$, as well. The heading angle is set to $\beta = 150$ degrees. The normalized forces can be seen in Figure \ref{fig:sphere_X_j} in terms of the real and imaginary parts of the two decomposed diffraction forces. In general, a reasonable level of agreement has been reached. Some discrepancies can be noted in the surge results for the real parts of the scattered solution towards the highest frequencies. This difference stems from the chosen pseudo-impulsive wave elevation, which should have contained a wider frequency span for the surge excitation force than the analogous heave case.

\subsection{A floating box}\label{sec:box}

The second benchmark study is of a floating box. The draft (submerged depth) of the box is $2$ [m], and the sides are both of length $2$ [m]; hereby, the characteristic length scale logically corresponds to $L = 2$ [m]. The water depth is set to $d = 3$ [m], and as for the sphere, we employ two symmetry planes with the same grading of the mesh. The polynomial order is set to $P = 4$. For the volumetric mesh, $N_{\text{elm}} = 1,084$ elements are used, where $132$ surfaces elements tessellate the body boundary of the quarter of the submerged box. With this, the system becomes of size $N_{\text{DOF}} = 17,035$.

\subsubsection{Added mass and damping coefficients}

\begin{figure}[t]
    \begin{minipage}[t]{.45\textwidth}
        \centering
        \includegraphics[width=0.85\textwidth]{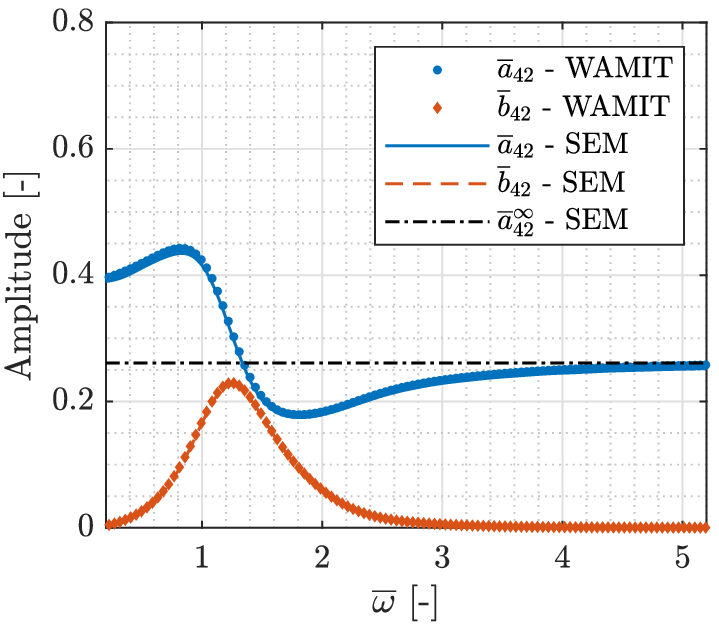}
        \subcaption{Roll-Sway}
        \vfill
        \centering
        \includegraphics[width=0.85\textwidth]{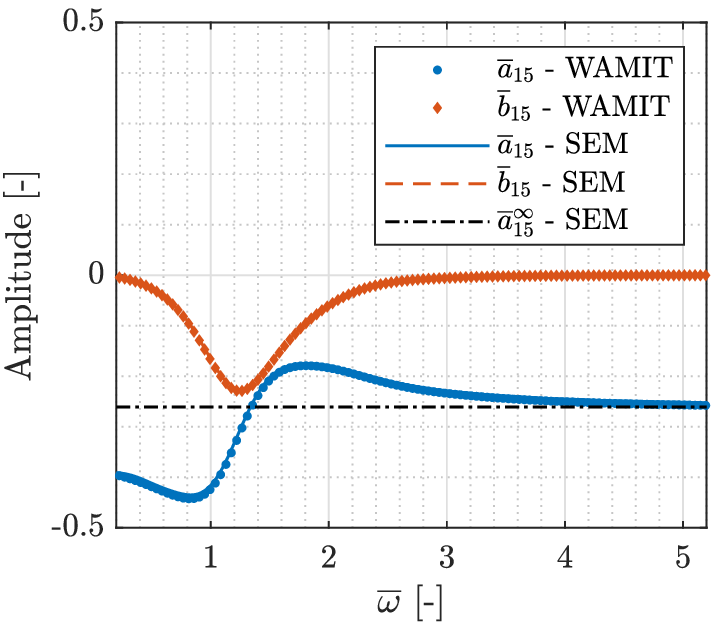}
        \subcaption{Surge-Pitch}
    \end{minipage}
    \hfill
    \begin{minipage}[t]{.45\textwidth}
        \centering
        \includegraphics[width=0.85\textwidth]{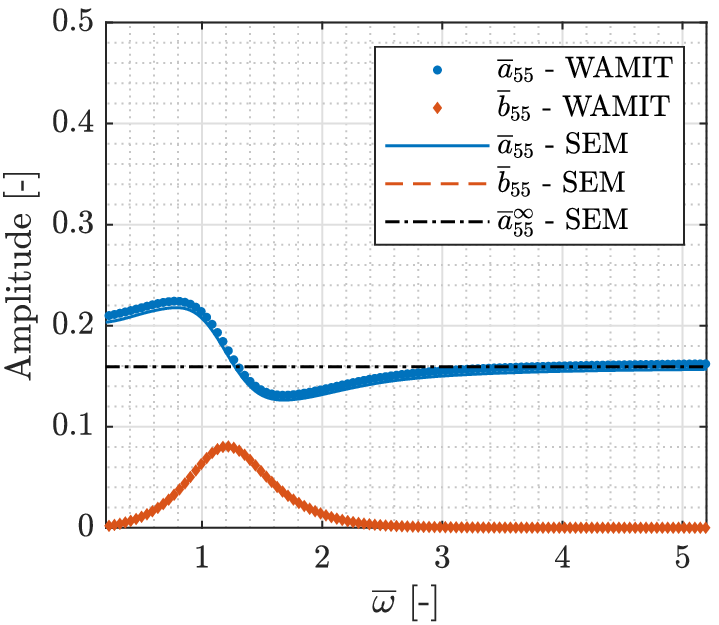}
        \subcaption{Pitch-Pitch}
        \vfill
        \centering
        \includegraphics[width=0.85\textwidth]{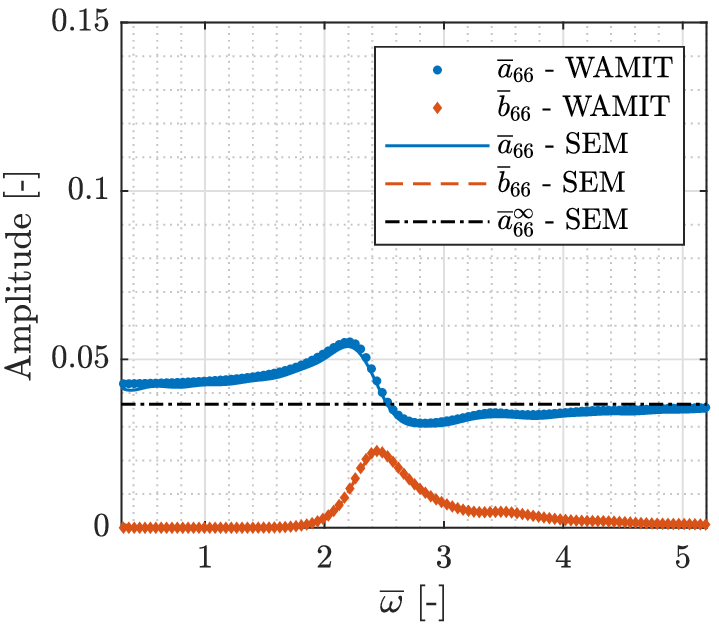}
        \subcaption{Yaw-Yaw}
    \end{minipage}  
    \caption{Added mass and damping coefficients, $\bar{a}_{jk}$ and $\bar{b}_{jk}$, for a floating box of length-scale $L = 2$ [m] on finite depth, $h = 3$ [m]. Comparison between the proposed SEM model and WAMIT results.}
    \label{fig:box_a_jk_b_jk}
\end{figure}

Now, solving the radiation problem with the pseudo-impulsive forcing in the sway, pitch, and yaw direction, and considering coupled and decoupled added mass and damping coefficients hereof is seen in Figure \ref{fig:box_a_jk_b_jk}. From the figure, the results match visually well with WAMIT. It can be argued that the added mass for the decoupled pitch motion, $\bar{a}_{55}$, is a bit under-resolved, hence slightly underestimating compared to the WAMIT results. The added mass clearly convergence towards the infinite added mass for all simulations.

\subsubsection{Wave excitation forces}

\begin{figure}[t]
    \begin{minipage}[t]{.45\textwidth}
        \centering
        \includegraphics[width=0.85\textwidth]{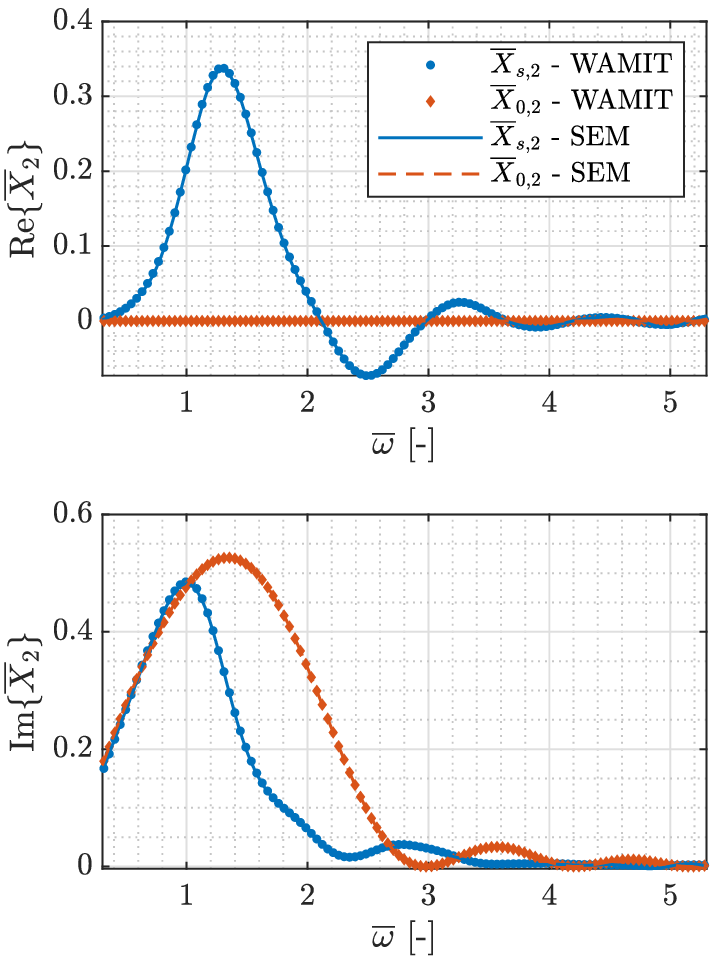}
        \subcaption{Sway}
    \end{minipage}
    \hfill
    \begin{minipage}[t]{.45\textwidth}
        \centering
        \includegraphics[width=0.85\textwidth]{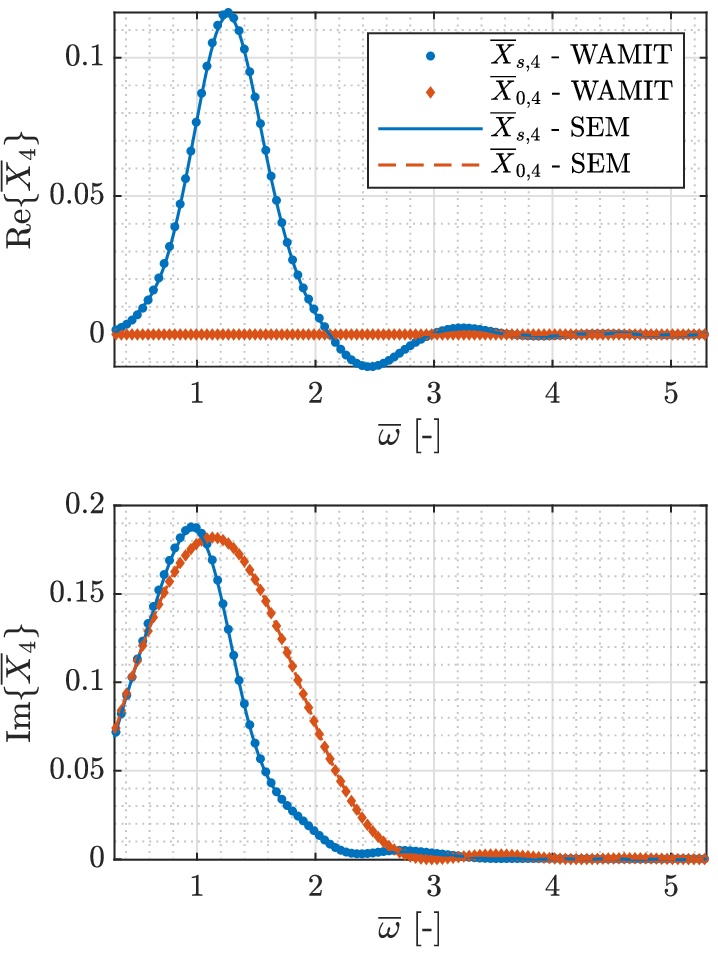}
        \subcaption{Roll}
    \end{minipage}  
    \caption{Wave excitation forces, $\bar{X}_{0,j}$ and $\bar{X}_{s,j}$ for a floating box of length-scale $L = 2$ [m] on finite depth, $h = 3$ [m], with heading angle $\beta = 135$ degrees. Comparison between the proposed SEM model and WAMIT results.}
    \label{fig:box_X_j}
\end{figure}

With the heading angle set to $\beta = 135$ degrees, the diffraction problem is solved, and highlighted results for sway and roll are shown in Figure \ref{fig:box_X_j}. Based on the figure, a satisfactory agreement can be observed when compared with the WAMIT results.

\subsection{An oscillating water column}\label{sec:OWC}

\begin{figure}[b] 
    \centering
    \includegraphics[width=0.45\textwidth]{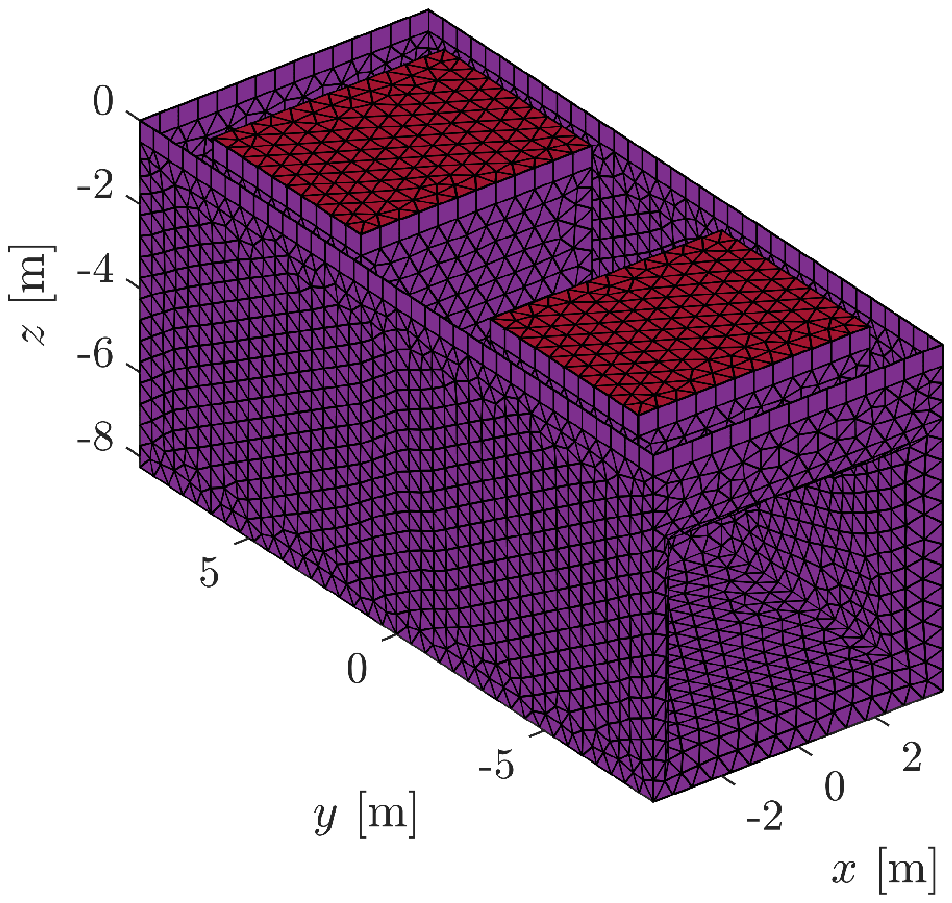}
    \caption{A surface mesh example of the OWC. Purple surface elements represent the body boundary, $\Gamma^{\text{body}}$, and the red elements represent the special internal body boundary, $\Gamma^{\text{s},1}$. Mesh is mirrored about $y = 0$ to show the complete OWC.}
    \label{fig:OWC_mesh}
\end{figure}

The last numerical experiment -- and perhaps most relevant in an applied context -- is of a one-way OWC-type wave energy converter as described and analyzed in \cite{joensen2023hydrodynamic}. The setup in the numerical wave tank is as follows: One symmetry plane is exploited about $y = 0$, and the length scale of the OWC is $L = 7.5$ [m]. Moreover, the draft is $8.25$ [m], the half-length (in the $y$-direction) is $8.7$ [m], and the width (in the $x$-direction) is $7.5$ [m]. The structure is placed on a water depth of $h = 32.5$ [m]. Inside the OWC, a chamber of width $6$ [m] and length $5$ [m] is located. In Figure \ref{fig:OWC_mesh}, a surface mesh example of the OWC is shown, where the purple surface elements represent the body boundary, $\Gamma^{\text{body}}$, and the red ones represent a special internal body boundary, $\Gamma^{\text{s},1}$ (more on this in the following). For the simulation, $N_{\text{elm}} = 19,338$ three-dimensional elements are considered, with $4,541$ surface elements representing the floating body, hereof $4,253$ on $\Gamma^{\text{body}}$ and $288$ on $\Gamma^{\text{s},1}$. This results in a system of size $N_{\text{DOF}} = 35,602$, when a polynomial order of $P = 2$ is used. 

\subsubsection{Generalized modes and special boundaries}

As announced in Section \ref{sec:mathematical_problem}, generalized modes and special boundaries are considered in the following. As for the latter, we employ one extra boundary, $\Gamma^{\text{s},1}$, which defines the internal surface of the OWC (red colored elements in Figure \ref{fig:OWC_mesh}). On this horizontal surface, we apply generalized modes according to \cite{newman1994wave}. Mathematically speaking, such modes are treated similarly to rigid body modes. For the OWC, $K = 2$ implies that two additional modes will be handled: a piston motion $(k = 7)$ and a first sloshing mode for the chamber $(k = 8)$. Each generalized mode is associated with a generalized normal, $n_k$, that reads
\begin{equation}
    n_7 = 1, \quad \text{and} \quad n_8 = \cos \left(\frac{\pi}{5} (y-2.2) \right),
\end{equation}
where $5$ [m] corresponds to the length of the chamber and $2.2$ [m] is the body-fixed position at the chamber's back wall. 

\subsubsection{Added mass and damping coefficients}

\begin{figure}[t]
    \begin{minipage}[t]{.45\textwidth}
        \centering
        \includegraphics[width=0.85\textwidth]{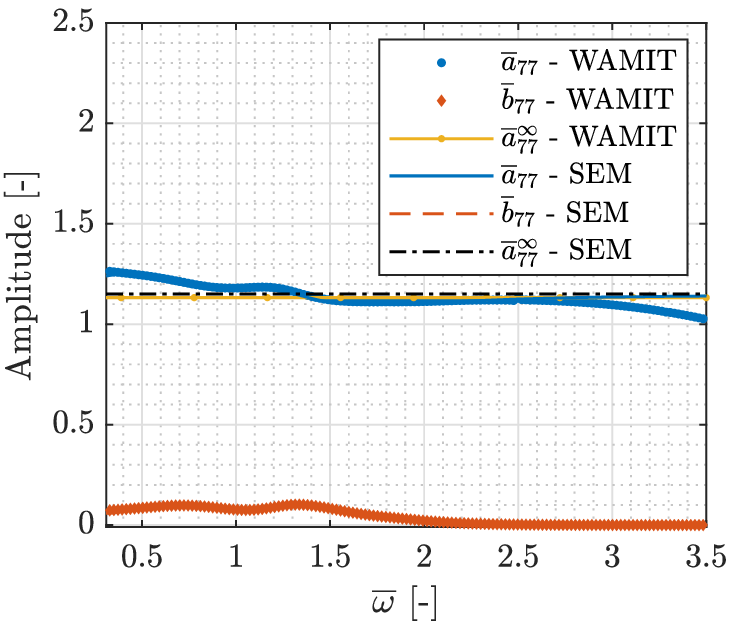}
        \subcaption{$(j,k) = (7,7)$}
    \end{minipage}
    \hfill
    \begin{minipage}[t]{.45\textwidth}
        \centering
        \includegraphics[width=0.85\textwidth]{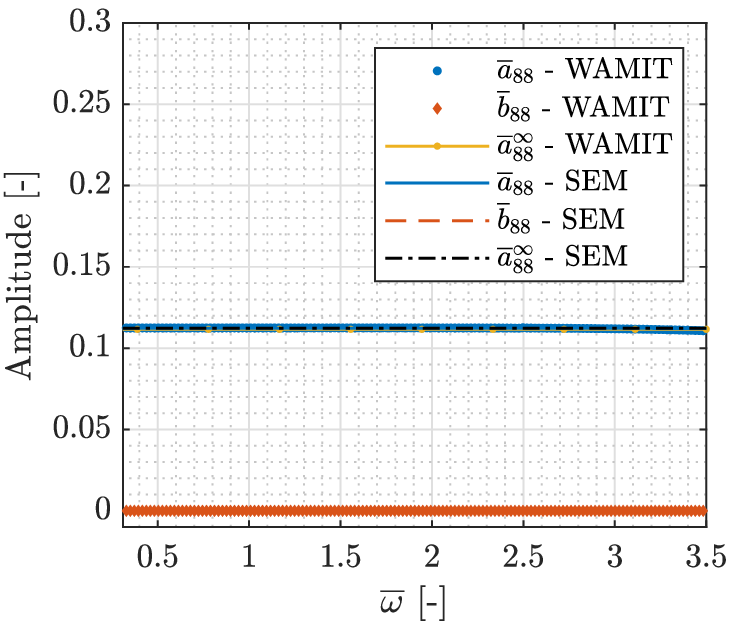}
        \subcaption{$(j,k) = (8,8)$}
    \end{minipage}  
    \caption{Added mass and damping coefficients, $\bar{a}_{jk}$ and $\bar{b}_{jk}$, for OWC-type structure. Comparison between the proposed SEM model and WAMIT results. For the piston motion $(k=7)$ and the first sloshing mode $(k=8)$.}
    \label{fig:OWC_a_jk_b_jk}
\end{figure}

For the radiation potential, we solve the problem twice -- one for each generalized mode -- where the results can be seen in Figure \ref{fig:OWC_a_jk_b_jk}. Here, it should be noted that the task of computing the infinite-frequency added mass, $a_{jk}^{\infty}$, for the OWC is a well-known singular problem when using WAMIT (or other BEM models). The singularity is due to ill-conditioning of the system when panels are placed in the plane of the free surface, as is the case with the internal surface, $\Gamma^{\text{s},1}$. The practical workaround, e.g., used in \cite{bingham2021ocean}, is to reflect the domain around $z=0$ and submerge it to $2h$. Another consequence of the singularity is related to the convergence of the added mass for large frequencies, which can be expected to be very slow. Luckily, this infinite-frequency problem is well-defined and trivial to solve using the SEM-based model.

From Figure \ref{fig:OWC_a_jk_b_jk}, we see the radiation coefficients for $(j,k) = (7,7)$ and $(j,k) = (8,8)$. The decoupled piston motion's coefficients are seen to agree well with the WAMIT simulation. However, for $\bar{\omega} > 2.5$, the SEM results differ due to the aforementioned singularity of the BEM solver, as the higher-frequency WAMIT added mass is not properly converged. Moreover, it can be observed that the decoupled first sloshing mode is \textit{waveless}, as the added mass is constant with zero damping. In both cases, the infinite-frequency added mass agrees reasonably.

 \begin{figure}[h]
    \begin{minipage}[t]{.45\textwidth}
        \centering
        \includegraphics[width=0.85\textwidth]{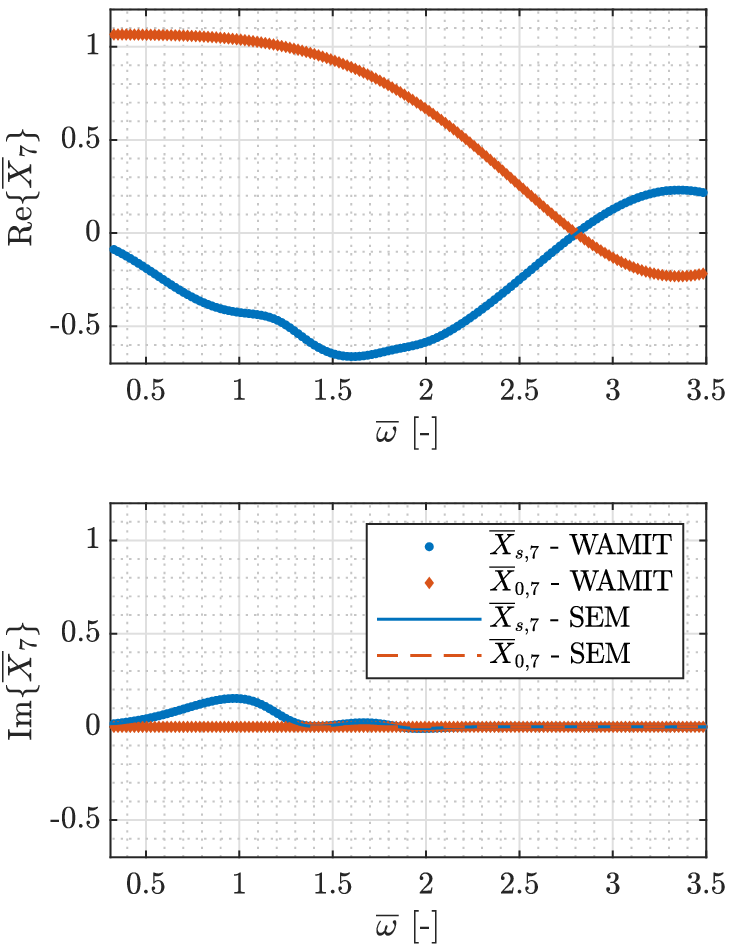}
        \subcaption{Piston mode.}
    \end{minipage}
    \hfill
    \begin{minipage}[t]{.45\textwidth}
        \centering
        \includegraphics[width=0.85\textwidth]{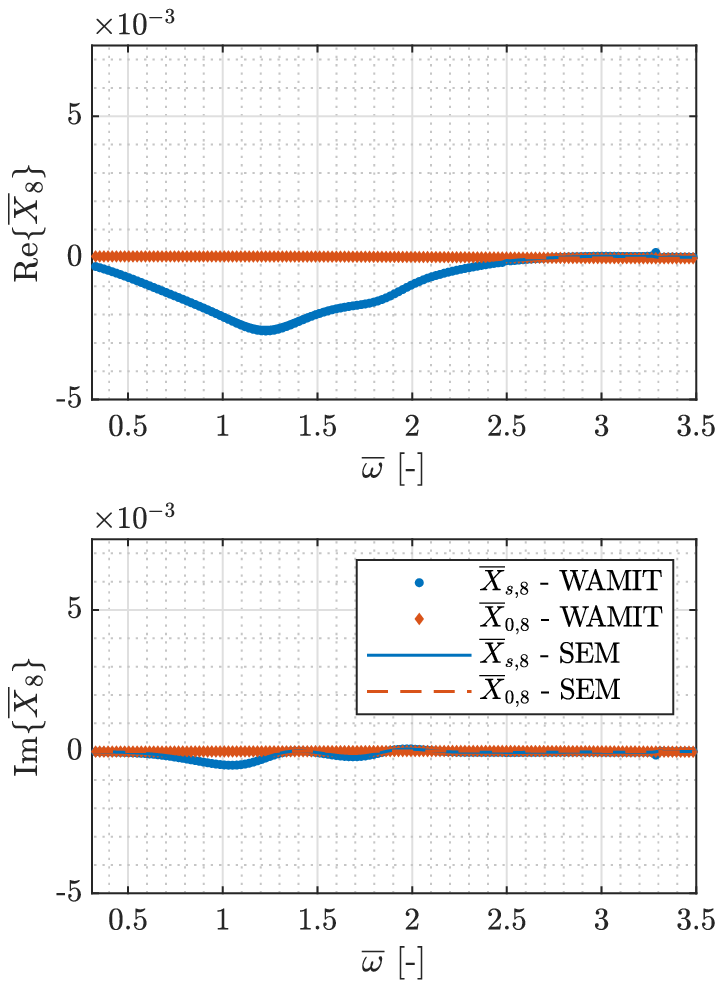}
        \subcaption{First sloshing mode.}
    \end{minipage}  
    \caption{Wave excitation forces, $\bar{X}_{0,j}$ and $\bar{X}_{s,j}$ for a for OWC-type structure with heading angle $\beta = \pi$. Comparison between the proposed SEM model and WAMIT results.}
    \label{fig:OWC_X_j}
\end{figure}

\subsubsection{Wave excitation forces}

For the diffraction problem, waves will be propagating along the $(y = 0)$ plane, hence $\beta = pi$. Moreover, only the symmetrical solution of the scattered potential will be used to compute the diffraction excitation forces. The results can be seen in Figure \ref{fig:OWC_X_j}, showing close to a perfect agreement with the BEM solver.

\section{Conclusion}\label{sec:conclusion}

This paper proposed a novel three-dimensional spectral element model for solving the linear radiation and diffraction problem using a pseudo-impulsive time-domain formulation. The governing equations were established for modeling floating offshore structures and their interaction with ocean waves approaching at an arbitrary heading angle, $\beta$. The high-order accurate spectral element solver is complemented by fourth-order accurate solutions for time integration and differentiation. Moreover, the FFT and its inverse were employed for the linear mapping between the time and frequency domains.

The expected convergence rate was confirmed for the discrete Laplace operator under mesh refinement and increase of the polynomial basis functions on \textit{hybrid} meshes. Lastly, various frequency-based quantities (added mass and damping coefficients and wave excitation forces) were simulated and compared with WAMIT results for a floating sphere, a floating box, and an OWC-type wave energy converter. The latter included generalized modes of a special internal surface to mimic the piston motion and the first sloshing mode. In general, the level of agreement was satisfactory in all cases, with little difference compared to the \textit{state-of-the-art} WAMIT solver. Moreover, the SEM model was naturally able to compute the infinite-frequency added mass coefficient for problems that are singular for WAMIT and other BEM-type solvers.

This paper constitutes the finalization of the complete linear radiation and diffraction problem using the spectral element method, where the two-dimensional model is outlined in \cite{visbech2023spectral} and initial three-dimensional results for the radiation problem can be seen in \cite{visbech2023high}. 

\subsection{Acknowledgments}

The research was conducted at the Technical University of Denmark (DTU) at the Department of Applied Mathematics and Computer Science in close collaboration with the Department of Civil and Mechanical Engineering (DTU Construct). The authors would like to express their gratitude to Ph.D. Mostafa Amini-Afshar from DTU Contruct for his help in discussing the boundary condition of the diffraction problem. Moreover, this work partly contributes to the activities of the research project: ”Efficient Added Mass Calculations for Large and Complex Floating Offshore Structures” supported by COWIfonden (Grant no. A-159.15). Lastly, the DTU Computing Center (DCC) provided access to computational resources.

\subsection{Declarations}

The authors of this paper hereby declare that there is no known personal, financial, or professional conflict of interest associated with completing the work presented in this paper.
\appendix

\section{Appendix}\label{sec:appendix}

\subsection{Appendix 1 - Body boundary conditions for the decomposed diffraction problem}\label{sec:appendix_1}

The key to understanding the decomposition in \eqref{eq:decomp} and how to construct the body boundary condition in \eqref{eq:body_bc_diffraction} for the diffraction problem when having decomposed it -- into symmetrical/anti-symmetrical and/or symmetrical/anti-symmetrical potentials -- is described in the following. 

First, considering the unit body normal vector, $\boldsymbol{n} = (n_x,n_y,n_z)$, and the symmetry and anti-symmetry relations hereof. Arguably
\begin{equation}
    \boldsymbol{n} = (\underbrace{n_x}_{\text{SA}},\underbrace{n_y}_{\text{AS}},\underbrace{n_z}_{\text{SS}}), \quad \text{on} \quad \Gamma^{\text{body}},
\end{equation}
using the same notation as introduced in Section \ref{sec:special}. Secondly, considering the gradient of the $i = \text{SS}$ scattered potential, $\nabla \phi_s^{\text{SS}}$, and symmetry and anti-symmetry relations hereof as
\begin{equation}
    \nabla \phi_s^{\text{SS}} =  (\underbrace{\partial_x \phi_s^{\text{SS}}}_{\text{AS}}, \underbrace{\partial_y \phi_s^{\text{SS}}}_{\text{SA}}, \underbrace{\partial_z \phi_s^{\text{SS}}}_{\text{SS}} ),
\end{equation}
and when combining $\partial_n  \nabla \phi_s^{\text{SS}} = \boldsymbol{n} \cdot \nabla \phi_s^{\text{SS}}$
\begin{equation}
     \boldsymbol{n} \cdot \nabla \phi_s^{\text{SS}} = \underbrace{n_x \partial_x \phi_s^{\text{SS}}}_{\text{SS}} + \underbrace{n_y \partial_y \phi_s^{\text{SS}}}_{\text{SS}} + \underbrace{n_z \partial_z \phi_s^{\text{SS}}}_{\text{SS}},  \quad \text{on} \quad \Gamma^{\text{body}},
\end{equation}
thus highlighting that the forcing of the SS part is double symmetrical as well. Such conditions can only be constructed when combining SS components from the incident wave, hereby yielding the body condition
\begin{equation}
    \partial_n \phi_s^i = - \partial_n \phi_0^i = - \boldsymbol{n} \cdot \IF{\nabla \Psi^i \F{\tilde{\zeta}_0}}, \quad \text{on} \quad \Gamma^{\text{body}}, \quad \text{for} \quad i = \{\text{SS, SA, AS, AA}\}.
\end{equation}

A similar exercise can be carried out for SA, AS, and SS, but is omitted in the presented paper.

\bibliographystyle{unsrtnat} 



\end{document}